\documentclass[12pt]{amsart}
\usepackage{amsmath,amssymb,amsfonts,amsthm,amsopn,dsfont}
\usepackage{graphics}
\usepackage{pb-diagram}

\newcommand{\ft}{Fourier transform}

\newcommand{\tf}{time-frequency}

\newcommand{\modsp}{modulation space}

\newtheorem{theorem}{Theorem}[section]
\newtheorem{corollary}[theorem]{Corollary}

\newtheorem{example}[theorem]{Example}
\newtheorem{lemma}[theorem]{Lemma}
\numberwithin{equation}{section}

\newtheorem{proposition}[theorem]{Proposition}
\newtheorem{remark}[theorem]{Remark}

\newcommand{\beqa}{\begin{eqnarray*}}
\newcommand{\eeqa}{\end{eqnarray*}}

\DeclareMathOperator*{\supp}{supp}

\newcommand{\field}[1]{\mathbb{#1}}
\newcommand{\bR}{\field{R}}        
\newcommand{\bZ}{\field{Z}}        
\def\omega{\eta}

\def\la{\lambda}

 \def\cF{\mathcal{F}}              
 \def\cS{\mathcal{S}}

 \def\cM{\mathcal{M}}

 \def\cC{\mathcal{C}}

\def\a{\aleph}

\def\rd{\bR^d}

\def\rdd{{\bR^{2d}}}

\def\lrd{L^2(\rd)}

\def\intrd{\int_{\rd}}

\def\R{\right)}
\def\l{\langle}
\def\r{\rangle}
\def\<{\left<}
\def\>{\right>}

\def\inv{^{-1}}

\def\mv1{M_v^1}

\def\Mmpq{M_m^{p,q}}
\def\phas{(x,\eta )}

\def\o{\omega}
\def\a{\alpha}

\def\R{\mathbb{R}}
\def\Ren{\mathbb{R}^d}
\def\Renn{\mathbb{R}^{2d}}

\def\sch{\mathcal{S}}

\def\Fur{\mathcal{F}}

\def\f{\varphi}

\def\Sn2{S_{2}(L^{2}(\Ren))}
\def\S1{S_{1}(L^{2}(\Ren))}
\def\sig00{\sigma_{0,0}}

\def\la{\langle}
\def\ra{\rangle}

\begin{document}

\title[Boundedness of Schr\"odinger type propagators on modulation spaces]{Boundedness of Schr\"odinger type propagators on modulation spaces}

\author{Elena Cordero and Fabio Nicola}
\address{Department of Mathematics,  University of Torino,
Via Carlo Alberto 10, 10123
Torino, Italy}
\address{Dipartimento di Matematica, Politecnico di
Torino, Corso Duca degli
Abruzzi 24, 10129 Torino,
Italy}
\email{elena.cordero@unito.it}
\email{fabio.nicola@polito.it}
\date{}

\begin{abstract}
It is known that Fourier
integral operators arising
when solving
Schr\"odinger-type operators
 are bounded
on the modulation spaces $\cM^{p,q}$, for $1\leq p=q\leq\infty$,
provided their symbols belong to the Sj\"ostrand class
$M^{\infty,1}$. However, they generally fail to be bounded on
$\cM^{p,q}$ for $p\not=q$. In this paper we study several
additional conditions, to be imposed on the phase or on the
symbol, which guarantee the boundedness on $\cM^{p,q}$ for
$p\not=q$, and between $\cM^{p,q}\to\cM^{q,p}$, $1\leq q<
p\leq\infty$. We also study similar problems for operators acting
on Wiener amalgam spaces, recapturing, in particular, some recent
results for metaplectic operators. Our arguments make heavily use
of the uncertainty principle.
\end{abstract}

\maketitle

\section{Introduction}
The paper is concerned with the study
of Fourier integral operators (FIOs)
defined by
\begin{equation}\label{fio}
    Tf(x)=\intrd e^{2\pi i \Phi\phas}
    \sigma\phas
    \hat{f}(\eta)d\eta,
\end{equation}
for, say, $f\in\mathcal{S}(\R^d)$. The
functions $\sigma$ and $\Phi$ are
called symbol and phase, respectively.
Here the \ft\, of $f$ is normalized to
be ${\hat
  {f}}(\eta)=\int
f(x)e^{-2\pi i x\eta}dx$. If
$\sigma\in L^\infty$ and the
phase $\Phi$ is real, the
integral converges absolutely
and defines a function in
$L^\infty$.\par The phase
function $\Phi(x,\eta)$
fulfills the following
properties:\\
(i) $\Phi\in \cC^{\infty}(\rdd)$;\\
(ii) there exist constants $C_\alpha>0$
such that
\begin{equation}\label{phasedecay}
|\partial^\a \Phi(x,\eta)|\leq
C_\a,\quad
\forall\a\in\bZ^{2d}_+,\,\,|\a|\geq
2;\end{equation}
 (iii)
 there exists $\delta>0$ such
that
\begin{equation}\label{detcond}
   \left|{\rm det}\,
\left(\frac{\partial^2\Phi}{\partial
x_i\partial \eta_l}\Big|_{
(x,\eta)}\right)\right|\geq
\delta\quad \forall
(x,\eta)\in \R^{2d}.
\end{equation}
Note that our phases differ from those
(positively homogeneous of degree 1 in
$\eta$) of FIOs arising in the solution
of hyperbolic equations (see, e.g.,
\cite{cordero-nicola-rodino,hormander71,ruzhsugimoto,seegersoggestein}).
Indeed, FIOs are a mathematical tool to
study a variety of problems in partial
differential equations, and our FIOs
arise naturally in the study of the
Cauchy problem for Schr\"odinger-type
operators (see, e.g.,
\cite{cordero2,fio1,numerico,folland,helffer84,helffer-rob1}).
Basic examples of phase functions
within the class under consideration
are quadratic forms in the variables
$x,\eta$ (see Example
\ref{metaplettici} below).
\par Continuing the study pursued in
\cite{fio1}, we focus on boundedness
results for these operators, when
acting on two classes of Banach spaces,
widely used in time-frequency analysis,
known as {\it modulation spaces}  and
{\it Wiener amalgam spaces}, denoted by
$M^{p,q}$ and $W(\Fur L^p,L^q)$,
respectively, with $1\leq p,q\leq
\infty$. To be definite, we recall the
definition of these spaces, introduced
by H. Feichtinger (see
\cite{F1,grochenig} and Section
\ref{preliminari} below for
details).\par
 In short, given a positive weight function $m$ on $\R^{2d}$, with $m\in\cS'(\rdd)$, we say
that a temperate distribution $f$ belongs to $M^{p,q}_m(\R^d)$,
$1\leq p,q\leq \infty$, if its short-time Fourier
 transform (STFT) $V_g
f\phas$, defined in \eqref{STFT} below, fulfills $V_g f\phas
m\phas\in L^{p,q}(\rdd)=L^q(\rd_{\o},L^p(\rd_x))$, with the norm
\begin{equation}\label{modulas}
\|f\|_{M^{p,q}_m}:=\|(V_gf)m\|_{L^q_{\o}L^p_x}<\infty.
\end{equation}
 Here $g$ is a non-zero
(so-called window) function in $\cS(\rd)$, which in \eqref{STFT}
is first translated and then multiplied by $f$ to localize $f$
near any point $x$. Changing $g\in \cS(\rd)$ produces equivalent
norms. Taking the two norms above in the converse order yields the
norm in $W(\Fur L^p,L^q)(\rd)$:
\begin{equation}\label{wiener}
\|f\|_{W(\Fur
L^p,L^q)}:=\|V_gf\|_{L^p_xL^q_{\o}}<\infty.
\end{equation}
 The spaces
$\mathcal{M}^{p,q}_m(\rd)$ and $\mathcal{W}(\Fur L^p,L^q)(\rd)$
are defined as
 the closure of $\cS(\rd)$ in the  ${M}^{p,q}_m(\rd)$ and $W(\Fur
L^p,L^q)$ norm, respectively. For heuristic purposes,
distributions in $M^{p,q}$, as well as in $W(\Fur L^q,L^p)$, may
be regarded as functions which are locally in ${\Fur L^q}$ and
decay at infinity like functions in $L^p$. Among their properties,
we highlight the important relation $W(\Fur L^p,L^q)=\Fur
(M^{p,q})$.\par The action on the spaces $\cM^p:=\cM^{p,p}$ of
FIOs as above already appeared in \cite{toft07,fio1} (see also
\cite{asada-fuji,benyi,boulkhemair,ruzhsugimoto}). It is  a basic
result that FIOs with symbols in the Sj\"ostrand class
$M^{\infty,1}(\rdd)$ extend to bounded operators on $\cM^p(\rd)$,
$1\leq p\leq\infty$. Applications to issues of classical analysis
were also given in \cite{cordero-nicola-rodino}. Moreover, it was
observed in \cite{fio1} that boundedness generally fails on the
spaces $\mathcal{M}^{p,q}(\rd)$, with $p\not=q$, although it can
hold under an additional condition on the phase function. The
present paper is devoted to a more systematic study of the
conditions which guarantee the boundedness on
$\mathcal{M}^{p,q}(\rd)$, for $p\not=q$.\par Our first result is
in fact a generalization of \cite[Theorem 11]{benyi} and
\cite[Theorem 5.2]{fio1}
 to
the case of rougher symbols.
\begin{theorem}\label{cont}
Consider a phase function $\Phi$ satisfying {\rm (i)}, {\rm (ii)},
and {\rm (iii)}, and a symbol $\sigma\in M^{\infty,1}(\rdd)$.
Suppose, in addition, that
\begin{equation}\label{fase}
\sup_{x,x',\eta\in\rd}\left|
\nabla_x\Phi(x,\eta)-\nabla_x\Phi(x',\eta)\right|<\infty.
\end{equation}
Then, the corresponding Fourier integral operator $T$ extends to a
bounded operator on $\mathcal{M}^{p,q}(\rd)$, for every $1\leq
p,q\leq\infty$.
\end{theorem}
The condition \eqref{fase} is seen to be essential for the
conclusion to hold. In fact it was shown in \cite[Proposition
7.1]{fio1} that the pointwise multiplication operator by
$e^{-\pi|x|^2}$ (which has phase
$\Phi(x,\eta)=x\eta-\frac{|x|^2}{2}$ and symbol $\sigma\equiv1$)
is not bounded on any $\mathcal{M}^{p,q}$, with $p\not=q$ (see
also Theorem \ref{3.1} below).\par If we drop the condition
\eqref{fase}, we need some further decay condition on the symbol,
as explained by the next result.\par
 For $s_1,s_2\in\R$, we define the weight function
 $v_{s_1,s_2}(x,\eta):=\langle
 x\rangle^{s_1}\langle\eta\rangle^{s_2}$, $(x,\eta)\in\R^{2d}$.
\begin{theorem}\label{main}
Consider a phase $\Phi$ satisfying {\rm (i)}, {\rm (ii)} and {\rm
(iii)}, and a symbol $\sigma\in
M_{v_{s_1,s_2}\otimes1}^{\infty,1}(\rdd)$,
$s_1,s_2\in\R$.\par\noindent
 $(i)$ Let $1\leq
p\leq\infty$. If $s_1,s_2\geq0$, $T$ extends to a bounded operator
on $\mathcal{M}^{p}(\rd)$.\par\noindent $(ii)$ Let $1\leq
q<p\leq\infty$. If $s_1>d\left(\frac{1}{q}-\frac{1}{p}\right)$,
$s_2\geq0$, $T$ extends to a bounded operator on
$\mathcal{M}^{p,q}(\rd)$.
\par\noindent $(iii)$  Let
$1\leq p<q\leq\infty$. If $s_1\geq0$,
$s_2>d\left(\frac{1}{p}-\frac{1}{q}\right)$, $T$ extends to a
bounded operator on $\mathcal{M}^{p,q}(\rd)$.\par
 In all cases,
\begin{equation}\label{normaoperatore}
\|Tf\|_{\mathcal{M}^{p,q}}\lesssim\|\sigma\|_{M_{v_{s_1,s_2}\otimes1}^{\infty,1}}\|f\|_{\mathcal{M}^{p,q}}.
\end{equation}
\end{theorem}
Although we do not exhibit a complete
set of counterexamples for the
thresholds arising in Theorem
\ref{main}, some examples are given in
Section \ref{schroedinger}, and show
that the thresholds are in fact the
expected ones (see Remark
\ref{aspettato}).\par

Results for boundedness of FIOs between weighted \modsp s are
attained as well (see Section \ref{continuitasumpqw}).\par
 We also turn  our
attention to the boundedness of FIOs as above from the modulation
space $\mathcal{M}^{p,q}(\rd)$ into $\mathcal{M}^{q,p}(\rd)$.
 This study is mostly suggested by
the special case of metaplectic operators (corresponding to a
quadratic phase and symbol $\sigma\equiv1$), which was
investigated in detail  in \cite{cordero2}; see also Example
\ref{metaplettici} below.
\begin{theorem}\label{corollario}
Let $1\leq q<p\leq\infty$. Consider a phase $\Phi(x,\eta)$
satisfying {\rm (i)}, {\rm (ii)} and {\rm (iii)}. Moreover, assume
one of the following conditions:\par (a) the symbol $\sigma\in
M^{\infty,1}(\rdd)$ and, for some $\delta>0$,
\begin{equation}\label{detcondbis0}
   \left|{\rm det}\,
\left(\frac{\partial^2\Phi}{\partial
x_i\partial x_l}\Big|_{
(x,\eta)}\right)\right|\geq
\delta\quad \forall
(x,\eta)\in \R^{2d},
\end{equation}\par (b) the symbol
$\sigma\in M^{\infty,1}_{v_{s,0}\otimes1}(\rdd)$, with
$s>d\left(\frac{1}{q}-\frac{1}{p}\right)$.\par\noindent
  Then, the
corresponding FIO $T$ extends to a bounded operator
$\mathcal{M}^{p,q}(\rd)\to \mathcal{M}^{q,p}(\rd)$.\par
\end{theorem}
The additional assumptions $(a)$ or
$(b)$ are essential to guarantee the
boundedness. A counterexample in this
connection is given in Proposition
\ref{controesempio} below. Moreover,
even under those conditions, $T$
generally fails to be bounded between
$\mathcal{M}^{p,q}\to
\mathcal{M}^{q,p}$ if $q>p$; see
Proposition \ref{controesempio0}
below.\par For the sake of brevity, in
this Introduction  we only established
our results for FIOs acting on
 modulation spaces. In the subsequent sections we shall provide
corresponding results for Wiener amalgam spaces (Corollaries
\ref{mainwiener} and \ref{corgen}).\par As in \cite{fio1}, the
proof of our results relies  on a formula expressing the Gabor
matrix of the FIO $T$  in terms of the STFT of its symbol $\sigma$
(see \eqref{symbstft} below). Here the novelty is provided by the
combination of this formula with  Schur-type tests and the
uncertainty principle, in the form of Bernstein's inequality and
some generalizations.

\vskip0.1truecm
 The paper is organized as follows. In Section
\ref{preliminari} we prove some preliminary results of classical
analysis and we also collect the basic definitions and properties
of modulation and Wiener amalgam spaces.  In Section
\ref{continuitasumpq} we recall from \cite{fio1} a useful formula
for the Gabor matrix of the operator $T$, and  use it to prove
Theorems \ref{cont} and \ref{main}. In Section
\ref{continuitasumpqw} we study the action of a FIO $T$ on
weighted modulation spaces. In Section \ref{continuitasumqp} we
prove Theorem \ref{corollario}. Finally, in Section
\ref{schroedinger} some examples related to the Schr\"odinger
operators are exhibited: they reveal to be useful tests for the
sharpness of the above results.
\par \vskip0.5truecm

\par
\textbf{Notation.} We define
$|x|^2=x\cdot x$, for
$x\in\Ren$, where $x\cdot
y=xy$ is the scalar product
on $\Ren$. The space of
smooth functions with compact
support is denoted by
$\cC_0^\infty(\rd)$, the
Schwartz class is
$\sch(\Ren)$, the space of
tempered distributions
$\sch'(\Ren)$.    The Fourier
transform is normalized to be
${\hat
  {f}}(\o)=\Fur f(\o)=\int
f(t)e^{-2\pi i t\o}dt$.
 Translation and modulation operators ({\it time and frequency shifts}) are defined, respectively, by
$$ T_xf(t)=f(t-x)\quad{\rm and}\quad M_{\o}f(t)= e^{2\pi i \o
 t}f(t).$$
We have the formulas $(T_xf)\hat{} =
M_{-x}{\hat {f}}$, $(M_{\o}f)\hat{}
=T_{\o}{\hat {f}}$, and
$M_{\o}T_x=e^{2\pi i x\o}T_xM_{\o}$.
The inner product of two functions
$f,g\in \lrd$ is $\la f, g \ra=\intrd
f(t) \overline{g(t)}\,dt$, and its
extension to $\cS'\times\cS$ will be
also denoted by  $\la \cdot, \cdot
\ra$. The notation $A\lesssim B$ means
$A\leq c B$ for a suitable constant
$c>0$, whereas $A \asymp B$ means
$c\inv A \leq B \leq c A$, for some
$c\geq 1$. The symbol $B_1
\hookrightarrow B_2$ denotes the
continuous embedding of the space $B_1$
into $B_2$. The open ball in $\R^d$ of
center $x$ and radius $R$ will be
denoted by $B(x,R)$.
\section{Preliminary
results}\label{preliminari}
\subsection{Bernstein inequalities}
The core of our proofs relies on the
classical Bernstein's inequality (see,
e.g., \cite{wolff}) and some of its
generalizations, described in what
follows. Recall that the ball of center
$x\in\rd$ and radius $R>0$ is denoted
by $B(x,R)$.
\begin{lemma}[Bernstein's inequality] Let $f\in\cS'(\rd)$ such that  $\hat{f}$ is supported in $B(0,R)$, and let  $1\leq p\leq q\leq\infty$. Then, there exists a positive constant $C$ (independent of $f$, $R$, $p,$ $q$), such that
\begin{equation}\label{bernstein}
\|f\|_q\leq C R^{d\left(\frac1p-\frac1q\right)}\|f\|_p.
\end{equation}
\end{lemma}
We shall use also the following generalizations
of the Bernstein's
inequalities.

\begin{lemma}\label{lemma}
Consider a mapping $v:\R^d\to\R^d$ satisfying $|v(x)|\lesssim
|x|$, $\forall x\in\R^d$. Then, for any $s\geq0$,
\[
\sup_{x\in\R^d}\{\langle x\rangle^{-s} |f(v(x))|\}\lesssim\int
\langle x\rangle^{-s} |f(x)|dx,
\]
for every function $f\in
\cS'(\rd)$, such that  $\hat{f}$ is supported in $B(\eta_0,1)$, for some $\eta_0\in\rd$.
\end{lemma}
\begin{proof} Take a Schwartz function
$g$, satisfying
$\hat{g}(\eta)=1$ for
$|\eta|\leq1$. If $\hat{f}$
is supported in $B(\eta_0,1)$
we have
$\hat{f}=\hat{f}T_{\eta_0}
\hat{g}$. Hence
\begin{equation}\label{pma}
\langle
x\rangle^{-s}|f(v(x))|\leq\int
\langle
x\rangle^{-s}|f(y)||g(v(x)-y)|dy.
\end{equation}
Now we have
\[
\langle y\rangle^{s}\lesssim \langle y-v(x)\rangle^{s}\langle
v(x)\rangle^s\lesssim \langle y-v(x)\rangle^{s}\langle x\rangle^s,
\]
so that
\[
\langle
x\rangle^{-s}|g(v(x)-y)|\leq
\langle x\rangle^{-s}\langle
v(x)-y\rangle^{-N}\lesssim
\langle y\rangle^{-s} \langle
v(x)-y\rangle^{-N+s}.
\]
Using this inequality in \eqref{pma}, with $N\geq s$, we attain
the desired conclusion.
\end{proof}
We recall from \cite[Proposition 5.5]{wolff} the following
localized version of  Bernstein's inequality.
\begin{lemma}\label{osservazione1} Let
$N\geq1$ be an integer and
$\varphi(x)=(1+|x|^2)^{-N}$.
For $R>0$, $x\in\R^{d}$, let
$\varphi_{x,R}(v)=\varphi(\frac{v-x}{R})$,
$v\in \R^d$. There exists
$C_N>0$ such that
\[
\sup_{y\in B(x,R)}|f(y)|\leq C_N\,
\mu(B(x,R)))^{-1}\intrd\varphi_{x,R}(v)|f(v)|\,dv,
\]
for every $f\in \cS'(\rd)$ such
that $\hat{f}$ is supported
in $B(\eta_0,1/R)$, for some $\eta_0\in\rd$,
where  $\mu$ is the
Lebesgue measure.\end{lemma}
\subsection{Schur-type tests} The next proposition collects  Schur-type
tests assuring the boundedness of  integral operators on the
mixed-norm spaces $L^q_\omega L^p_x:=L^q(\R^d_\o;L^p(\R^d_x))$,
$1\leq p,q\leq \infty$.

\begin{proposition}\label{proschur}
Consider an integral operator $A$ on $\R^{2d}$, given by
\[
(A f)(x',\omega')=\iint K(x',\o';x,\o)f(x,\o)\,dx\,d\omega.
\]
(i) If $K\in L^\infty_{\o}L^1_{\o'}L^\infty_{x'}L^1_{x}$, then $A$
is continuous on
$L^1_{\o}L^\infty_{x}$.\\
(ii) If $K\in L^\infty_{\o'}L^1_{\o}L^\infty_{x}L^1_{x'}$, then
$A$ is continuous on
$L^\infty_{\o}L^1_{x}$.\\
(iii) If $K\in L^\infty_{\o}L^1_{\o'}L^\infty_{x'}L^1_{x}\cap
L^\infty_{\o'}L^1_{\o}L^\infty_{x}L^1_{x'}$, and, moreover, $K\in
L^\infty_{x',\o'}L^1_{x,\o}\cap L^\infty_{x,\o}L^1_{x',\o'}$, then
the operator $A$ is continuous on $L^q_{\o}L^p_{x}$, for
every $1\leq p,q\leq\infty$.\\
(iv) If $K\in L^\infty_{\o,\o'}L^1_{x,x'}$, then $A$ is continuous
$L^1_{\o}L^\infty_x\to L^{\infty}_{\o}L^1_{x}$.
\end{proposition}
\begin{proof}
 The
proof of all items, but (iv),
is just a repetition, with
obvious changes, of that of
\cite[Proposition 5.1]{fio1},
where a discrete version was
presented.\\ Let us now prove
(iv). We have
\begin{align*}
\|A
f\|_{L^\infty_{\o'}L^1_{x'}}&=\sup_{\o'\in\R^{d}}\int_{\R^d_{x'}}\left|\iint
K(x',\o';x,\o)f(x,\o)\,dx\,d\omega\right|\,dx'\\
&\leq \sup_{\o'\in\R^{d}}
\iiint_{\R^{3d}_{x,x',\o}}|K(x',\o';x,\o)f(x,\o)|
dx\,dx'\,d\o\\
&\leq \sup_{\o'\in\R^{d}}\int_{R^d_{\o}}\sup_{x\in\R^d}
|f(x,\o)|\iint|K(x',\o';x,\o)|dx\,dx'\,d\o\\
&\leq \|K\|_{L^\infty_{\o,\o'}L^1_{x,x'}}
\|f\|_{L^1_{\o}L^\infty_{x}}.
\end{align*}
This concludes the proof.
\end{proof}\\

\subsection{Modulation
spaces}\cite{F1,fe89-1,feichtinger90,
Feichtinger-grochenig89,grochenig,triebel83}
For $s\in\R$, we denote by
$\la\cdot\ra^s=(1+|\cdot|^2)^{s/2}$. In
what follows we limit ourselves to the
class of weight functions
$v_{s_1,s_2}\phas=\la x\ra^{s_1}\la
\o\ra^{s_2}$, $s_i\in\R$, $i=1,2$, on
$\rdd$, or $m(x,\eta,z,\zeta)=\la
x\ra^{s_1}\la
\o\ra^{s_2}=(v_{s_1,s_2}\otimes
1)(x,\eta,z,\zeta)$, on $\R^{4d}$. In
order to define such spaces, we make
use of the following \tf
\,representation: the short-time
Fourier transform (STFT) $V_gf$ of a
function/tempered distribution
$f\in\cS'(\rd)$ with respect to the the
window $g\in\cS(\rd)$ is defined by
\begin{equation}\label{STFT}
V_g f(x,\o)=\la f,M_\o T_xg\ra= \intrd e^{-2\pi i \o
y}f(y)\overline{g(y-x)}\,dy,
\end{equation}
i.e.,  the  Fourier transform $\cF$ applied to $f\overline{T_xg}$.

Given a non-zero window
$g\in\sch(\Ren)$, a weight function $m$
as those quoted above, and $1\leq
p,q\leq \infty$, the {\it
  modulation space} $M^{p,q}_m(\Ren)$ consists of all tempered
distributions $f\in\sch'(\Ren)$ such that $V_gf\in L^{p,q}_m(\Renn
)$ (weighted mixed-norm spaces). The norm on $M^{p,q}_m$ is
$$
\|f\|_{M^{p,q}_m}=\|V_gf\|_{L^{p,q}_m}=\left(\int_{\Ren}
  \left(\int_{\Ren}|V_gf(x,\o)|^pm(x,\o)^p\,
    dx\right)^{q/p}d\o\right)^{1/p}  \,
$$
(with obvious changes when $p=\infty$ or $q=\infty$). If $p=q$, we
write $M^p_m$ instead of $M^{p,p}_m$, and if $m(z)\equiv 1$ on
$\Renn$, then we write $M^{p,q}$ and $M^p$ for $M^{p,q}_m$ and
$M^{p,p}_m$, respectively. Then  $\Mmpq (\Ren )$ is a Banach space
whose definition is independent of the choice of the window $g$.
For the properties of these spaces we refer to the literature
quoted at the beginning of this subsection.

We define by $\mathcal{M}^{p,q}_m(\rd)$
 the closure of $\cS(\rd)$ in the  $M^{p,q}_m$-norm. Observe
 that $\mathcal{M}^{p,q}_m=M^{p,q}_m$, whenever the indices $p$ and
 $q$ are finite. They enjoy the duality property:
 $(\mathcal{M}^{p,q}_m)^*=\mathcal{M}^{p',q'}_{1/m}$, with $1<
 p,q<\infty$, and $p',q'$ being the conjugate exponents.

\par We  recall the inversion formula for
the STFT (see e.g. (\cite[Corollary 3.2.3]{grochenig}): if
$\|g\|_{L^2}=1$ and, for example, $u\in{L^2}(\R^d)$, it turns out
\begin{equation}\label{treduetre}
u=\int_{\R^{2d}} V_g u(x,\o)M_\o T_x g\, dx\,d\o.
\end{equation}

The following inequality, proved in \cite[Lemma
11.3.3]{grochenig},  is useful for changing windows.
\begin{lemma}\label{changewind}
Let
$g_0,g_1,\gamma\in\cS(\rd)$
such that $\la \gamma,
g_1\ra\not=0$ and let
$f\in\cS'(\rd)$. Then,
$$|V_{g_0} f\phas|\leq\frac1{|\la\gamma,g_1\ra|}(|V_{g_1} f|\ast|V_{g_0}\gamma|)\phas,
$$
for all $\phas\in\rdd$.
\end{lemma}
The complex interpolation theory for these spaces reads as follows
(see, e.g., \cite{Feichtinger-grochenig89}).
\begin{proposition}\label{interp} Let $0<\theta<1$,
 $p_j,q_j\in [1,\infty]$ and $m_j$ be $v$-moderate
 weight functions (i.e., $m_j(x+y)\leq C v(x)m_j(y)$, $v$ being
a submultiplicative weight), $j=1,2$.
Set
$$ \frac1p=\frac{1-\theta}{p_1}+\frac{\theta}{p_2},\quad \frac1q=\frac{1-\theta}{q_1}+\frac{\theta}{q_2},\quad m=m_1^{1-\theta}m_2^{\theta},$$
then
$$(\mathcal{M}^{p_1,q_1}_{m_1}(\rd),\mathcal{M}^{p_2,q_2}_{m_2}(\rd))_{[\theta]}=\mathcal{M}^{p,q}_{m}(\rd).$$
\end{proposition}

\begin{remark}\label{2.3}\rm
We observe that our results
are established as the
existence of a bounded
extension
$\mathcal{M}^{p,q}\to\mathcal{M}^{\tilde{p},\tilde{q}}$
of a class of FIOs $T$ with
symbols $\sigma$ in a
weighted space $M^{\infty,1}_m$, $m\geq 1$. It is
important to observe that
such an extension follows
from a uniform estimate of
the type
\begin{equation}\label{tipo}
\|Tf\|_{{M}^{\tilde{p},\tilde{q}}}\leq
C\|\sigma\|_{M^{\infty,1}_m}\|f\|_{{M}^{p,q}},\qquad
\forall f\in \cS(\R^d).
\end{equation}
Indeed, this estimate shows that $T$ extends to a bounded operator
$\mathcal{M}^{p,q}\to {M}^{\tilde{p},\tilde{q}}$. In order to
prove that this extension takes values in $\mathcal{M}^{\tilde{p},\tilde{q}}$, it suffices to verify that
$Tf\in\mathcal{M}^{\tilde{p},\tilde{q}}$ when $f$ is a Schwartz
function. This follows from \cite[Theorem 6.1]{fio1}. Indeed, if $f\in M^1$, then $Tf\in M^1$. 
\par Hence in the subsequent
proofs we will prove
estimates of the type
\eqref{tipo}.
\end{remark}
Boundedness
results dealing with FIOs having
symbols in weighted
modulation spaces and acting
on unweighted modulation
spaces could be
rephrased as
boundedness results for FIOs
with symbols in unweighted
spaces and acting on weighted
spaces, as explained below.
\begin{proposition}\label{proposizione1}
Let  $T$ be a FIO with symbol $\sigma$ and $\tilde{T}$ a FIO with the same
phase as $T$ and symbol
$$\tilde{\sigma}(x,\eta):=\l
x\r^{s_1}\sigma(x,\eta)\l\eta\r^{s_2},\quad x,\o \in\rd,\,\, s_i\in\R,\,\,i=1,2.$$
 Then,\\
(i) the operator $T$ is bounded from $\cM^{p,q}$ into
$\cM^{\tilde{p},\tilde{q}}$ if and only if the operator
$\tilde{T}$ is bounded from
$\cM^{p,q}_{v_{0,s_2}}$ into $\cM^{\tilde{p},\tilde{q}}_{v_{-s_1,0}}$.\\
(ii) It holds true
\begin{equation}\label{pesi}
\sigma\in
M^{\infty,1}_{v_{s_1,s_2}\otimes1}(\R^{2d})\Longleftrightarrow
\tilde{\sigma}\in
M^{\infty,1}(\R^{2d}).
\end{equation}
\end{proposition}
\begin{proof}
The proof is an immediate
consequence of \cite[Theorem
2.2, Corollary
2.3]{Toftweight}.  Indeed,
they guarantee that the
vertical arrows of the
following commutative diagram
define isomorphisms:
\begin{equation}\label{diagramma}
\begin{diagram}
\node{{\cM}^{p,q}}
\arrow{c,t}{T} \node{
{\cM}^{\tilde{p},\tilde{q}}}\arrow{s,r}{\text{pointwise
product by}\ \l x\r^{s_1}}
  \\
\node{{\cM}^{p,q}_{v_{0,s_2}}}\arrow{n,l}{\l
D\r^{s_2}}
\arrow{c,t}{\tilde{T}} \node{
{\cM}^{\tilde{p},\tilde{q}}_{v_{-s_1,0}}}.
\end{diagram}
\end{equation}
Moreover, the product by the weight
function
$\o_0(x,\eta,\zeta_1,\zeta_2):=\la
x\r^{s_1}\la\eta\r^{s_2}$, $x,\eta,$
$\zeta_1,\zeta_2\in\rd$, defined on
$\R^{4d}$, is a isomorphism from
$M^{\infty,1}_{v_{s_1,s_2}\otimes1}(\R^{2d})$
to $M^{\infty,1}(\R^{2d})$, that is
$(ii)$.
\end{proof}

\subsection{Wiener amalgam
spaces}\cite{F1,feichtinger90,fournier-stewart85} For $1\leq
p\leq\infty,$ $s\in\R$, we denote by
$L^p_s(\rd):=L^p_{\la\cdot\ra^s}(\rd)$. Let $ \cF L^p_s(\rd)$ be
the space of $f\in\cS'(\rd)$ such that $\hat{f}\in L^p_s(\rd)$.
Let $g \in \cS (\rd )$ be a non-zero window function. For $1\leq
p,q,\leq\infty$, the {\it Wiener amalgam space} $W(\cF
L^p_{s_1},L^q_{s_2})(\rd)$ with local component $\cF L^p_{s_1}$
and global component  $L^q_{s_2}$, $s_i\in\R$, $i=1,2$,  is
defined as the space of all functions or distributions for which
the norm
$$
\|f\|_{W(\cF L^p_{s_1},L^q_{s_2})}=\left( \int _{\Ren } \left(
\int _{\Ren } |\cF ({f\cdot T_{x}g})(\eta)|^p\la\eta\ra^{p s_1}
d\eta \right)^{q/p} \la x\ra^{qs_2}\, dx \right)^{1/q}
$$
is finite. Analogously to modulation spaces, we define by
${\mathcal {W}}(\cF L^p_{s_1},L^q_{s_2})(\rd)$ the closure of
$\cS(\rd)$ with respect to the $W(\cF L^p_{s_1},L^q_{s_2})$-norm.

 The  properties of Wiener amalgam
spaces are similar to those of  \modsp s, since we have the
following relation:
\begin{equation}\label{wienermod}
\cF (\mathcal{M}^{p,q}_{v_{s_1,s_2}})={\mathcal {W}}(\cF
L^p_{s_1},L^q_{s_2}),\quad s_i\in\R,\,i=1,2.
\end{equation}

We recall
from \cite[Lemma 5.3]{CN3}
the following auxiliary
result.
\begin{lemma}\label{lemm2}
For $a,b\in\R$, $a>0$, set
$G_{a+ib}(x)=(a+ib)^{-d/2}e^{-\frac{\pi|x|^2}{a+ib}}$,
and choose the Gaussian
$g(y)=e^{-\pi |y|^2}$ as
window function. Then, for
every $1\leq q,r\leq\infty$,
it follows that
\begin{equation}\label{lem1000}
|\widehat{G_{a+ib} T_x
g}(\omega)|=((a+1)^2+b^2)^{-d/4}
e^{-\frac{\pi}{(a+1)^2+b^2}[(a(a+1)+b^2)|
\omega|^2+2bx\omega+(a+1)|x|^2]},
\end{equation}
\begin{equation}\label{lem10}
\|G_{a+ib} T_x g\|_{\Fur
L^q}=\frac{\left((a+1)^2+b^2\right)^{\frac{d}{2}\left(\frac{1}{q}-\frac{1}{2}\right)}}{q^{\frac{d}{2q}}\left(a(a+1)+b^2\right)^{\frac{d}{2q}}}
e^{-\frac{\pi
a|x|^2}{a(a+1)+b^2}},
 \end{equation}
 and
\begin{equation}\label{lem2}
\|G_{a+ib}\|_{W(\Fur
L^q,L^r)}\asymp
\frac{\left((a+1)^2+b^2\right)^{
\frac{d}{2}\left(\frac{1}{q}-\frac{1}{2}
\right)}}{a^\frac{d}{2r}\left(a(a+1)+b^2\right)^{\frac{d}{2}\left(\frac{1}{q}-\frac{1}{r}\right)}}.
\end{equation}
\end{lemma}

\begin{lemma}\label{lemm3} Let $\sigma\in
M^{p,q}_{v_{s_1,s_2}\otimes 1}(\rdd)$, $s_i\in\R$, $i=1,2$, and
consider its adjoint $\sigma^*$ defined by
\begin{equation}\label{adj}
\sigma^\ast(x,\eta)=\overline{\sigma(\eta,x)}, \quad (x,\eta)\in\rdd.
\end{equation}
Then, $\sigma^*\in M^{p,q}_{v_{s_2,s_1}\otimes 1}(\rdd)$.
\end{lemma}
\begin{proof}
Fix a window function $g\in\cS(\rd)$. By a direct
computation,
\begin{align*}
V_g\sigma^\ast(z_1,z_2;\zeta_1,\zeta_2)&=
\int e^{-2\pi i(
y_1\zeta_1+y_2\zeta_2)}\overline{\sigma(y_2,y_1)}
\,\overline{g(y_1-z_1,y_2-z_2)}\,dy_1\,dy_2\\
&=\int e^{-2\pi i(
v_2\zeta_1+v_1\zeta_2)}\overline{\sigma(v_1,v_2)}
\,\overline{g(v_2-z_1,v_1-\zeta_2)}\,dv_1\,dv_2\\
&=V_{^tg}\overline{\sigma}(z_2,z_1;\zeta_2,\zeta_1)
\end{align*}
where we used the notation $^tg(y_1,y_2):=g(y_2,y_1)$. Since
$\overline{\sigma}\in M^{p,q}_{v_{s_1,s_2}\otimes 1}(\rdd)$, the
result immediately follows by the indipendence of the window
function for the computation of the modulation space norm.
\end{proof}

\section{Boundedness of
FIOs on $\mathcal{M}^{p,q}$}
 \label{continuitasumpq}

 In the
following we assume that the
phase function $\Phi$
satisfies the assumptions
(i), (ii) and (iii) in the
Introduction, so that we shall repeatedly use the following property.
\begin{remark}\label{cambiovariabili}\rm
It follows from the
assumptions (i), (ii) and
(iii) and Hadamard's global
inversion function theorem
(see, e.g., \cite{krantz})
that the mappings
\[
x\longmapsto\nabla_\eta\Phi(x,\eta)
\]
and
\[
\eta\longmapsto\nabla_x\Phi(x,\eta)
\]
are global diffeomorphisms of
$\R^d$, and their inverse
Jacobian determinant are
uniformly bounded with
respect to all variables.
\end{remark}

The proof of our results relies
 on a formula,  obtained in \cite[Section 6]{fio1}, which
expresses the Gabor matrix of
the FIO $T$  in
terms of the STFT of its
symbol $\sigma$. Namely,
choose a non-zero window
function $g\in \cS(\R^d)$,
and define
\begin{equation}\label{window}
\Psi_{(x',\o)}(y,\zeta):=e^{2\pi
i\Phi_{2,(x',\o)}(y,\zeta)}(\bar{g}\otimes\hat{g})(y,\zeta),\quad
(y,\zeta)\in\rdd
\end{equation}
with
\begin{equation}
  \label{eq:c11}
\Phi_{2,(x',\o)}(y,\zeta)=2\sum_{|\a|=2}\int_0^1(1-t)\partial^\a
\Phi((x',\o)+t(y,\zeta))\,dt\frac{(y,\zeta)^\a}{\a!}.
\end{equation}
Moreover, let
$$g_{x,\omega}(t)=\left(M_\omega
T_x g\right)(t),\quad t,x,\omega\in\R^d.$$
\begin{proposition}It turns out
\begin{equation}\label{symbstft}
 |\la T g_{x,\o},
g_{x',\o'}\ra| = |V_{\Psi_{(x',\o)}}
\sigma((x',\o),(\o'-\nabla_x\Phi(x',\o),x-\nabla_\eta
\Phi(x',\o)))|.
\end{equation}
\end{proposition}

\begin{remark}\rm\label{osser1}
Observe that the window $\Psi_{(x',\o)}$ of the STFT above depends
on the pair $(x',\o)$. However, the assumptions \eqref{phasedecay}
imply that these windows belong to a bounded subset of the
Schwartz space, i.e. the corresponding seminorms are uniformly
bounded with respect to $(x',\o)$.
\end{remark}

In this section we focus on the proofs
 of Theorems \ref{cont} and \ref{main}.
We first prove Theorem
\ref{cont}. We need the
following auxiliary results.

\begin{lemma}\label{changewind2}
Let $\Psi_0\in\cS(\rdd)$ with $\|\Psi_0\|_{L^2}=1$ and
$\Psi_{(x',\o)}$ be defined by \eqref{window}, $(x',\o)\in
\R^{2d}$, and $g\in\cS(\rd)$. Then,
\begin{equation}\label{newwind}
\int_{\bR^{4d}}\sup
_{(x',\o)\in\R^{2d}}|V_{\Psi_{(x',\eta)}}\Psi_0(w)|\,dw<\infty.
\end{equation}
\end{lemma}
\begin{proof} We shall
show that
\begin{equation}\label{lw} |V_{\Psi_{(x',\o)}}\Psi_0(w)|\leq C\la w\ra^{-(4d+1)},\quad \quad \forall (x',\o)\in\R^{2d}.
\end{equation}
Using the switching property of the STFT: $$
(V_fg)(x,\eta)=e^{-2\pi i\eta x}\overline{(V_gf)(-x,-\eta)},$$ we
observe that $|V_{\Psi_{(x',\eta)}}\Psi_0(w_1,w_2)|=
|V_{\Psi_0}\Psi_{(x',\o)}(-w_1,-w_2)|$, and by the even property
of of the weight $\la \cdot\ra$, relation \eqref{lw} is equivalent
to
\begin{equation}\label{l1} |V_{\Psi_0}\Psi_{(x',\eta)}(w)|\leq C\la w\ra^{-(4d+1)},\quad \quad \forall (x',\o)\in\rdd.
\end{equation}
Now, the mapping $V_{\Psi_0}$ is
continuous from $\cS(\rdd)$ to
$\cS(\bR^{4d})$ (see \cite[Chap.
11]{grochenig}), which combined with
the Remark \eqref{osser1} yields
\eqref{l1}.
\end{proof}\\
\begin{lemma}\label{osservazione2}
With the notation of
Lemma \ref{osservazione1}, for
every $R>0$ there exists
$C_R$ such that
\[
\iint\varphi_{x,R}(v)f(v)dv\,dx=C_R\int
f(x)\,dx,
\]
for every measurable function
$f\geq0$.
\end{lemma}
\begin{proof}
The proof is just an
application of Fubini's
theorem on the exchange of
integrals, since
$\int\varphi_{x,R}(v)dx=C_R$
is independent of $v$.
\end{proof}
\begin{proposition}\label{osservazione3}
Let $\Psi_0\in\cS(\R^{2d})$
supported in the ball
$B(0,1/R)$. Then
\[
\sup_{v\in
B(0,R)}|V_{\Psi_0}\sigma(u;z_1+v,z_2)|\leq
C_N \,\mu
(B(z_1,R))^{-1}\int\varphi_{z_1,R}(v)|V_{\Phi_0}\sigma(u;v,z_2)|dv,
\]
for every $u\in\R^{2d}$,
$z_1,z_2\in\R^d$.
\end{proposition}
\begin{proof}
It suffices to apply
Proposition
\ref{osservazione1} to the
function $\R^d\ni
v\longmapsto
V_{\Psi_0}\sigma(u;v,z_2)$.
Indeed, setting
$u=(u_1,u_2)\in\R^d\times\R^d$,
its Fourier transform has
support contained in the ball
$B(u_1,1/R)$.
\end{proof}
\begin{proof}[Proof of
Theorem \ref{cont}] Since the boundedness on
$\mathcal{M}^p=\mathcal{M}^{p,p}$, $1\leq p\leq\infty$, was
already proved in \cite[Theorem 6.1]{fio1}, see also \cite[Theorem
2.1]{toft07}, the desired result will follow by interpolation from
the cases $(p,q)=(\infty,1)$ and $(p,q)=(1,\infty)$. To prove
them, we observe that the inversion formula \eqref{treduetre} for
the STFT gives
\[
V_g (Tu)(x',\omega')=\int_{\R^{2d}}\langle Tg_{x,\omega},
g_{x',\o'}\rangle V_g u(x,\omega)dx\,d\omega.
\]
By Remark \ref{2.3}, the desired estimate therefore follows if we
prove that the map $K_T$ defined by
\[
K_T G(x',\omega')=\int_{\R^{2d}}\langle Tg_{x,\omega},
g_{x',\o'}\rangle G(x,\omega)dx\,d\omega
\]
is continuous on $L^1_\omega L^\infty_x$ and $L^\infty_\omega
L^1_x$. By Proposition \ref{proschur} it suffices to prove that
its integral kernel
\[
K_T(x',\omega';x,\omega)=\langle Tg_{x,\omega}, g_{x',\o'}\rangle
\]
satisfies
\begin{equation}\label{schur1}
K_T\in L^\infty_{\o}L^1_{\o'}L^\infty_{x'}L^1_{x},
\end{equation}
and
\begin{equation}\label{schur2}
K_T\in L^\infty_{\o'}L^1_{\o}L^\infty_{x}L^1_{x'}.
\end{equation}
Let us verify \eqref{schur1}.
By \eqref{symbstft} we have
\[
|K_T(x',\o';x,\o)|=|V_{\Psi_{(x',\o)}}\sigma(z(x',\o',x,\o))|,
\]
where \[
z(x',\o',x,\o)=(x',\o,\o'-\nabla_x\Phi(x',\o),x-\nabla_\eta
\Phi(x',\o)). \]
 By Lemma
\ref{changewind} for
$g_1=\gamma=\Psi_0$, $\Psi_0\in \cS(\rdd)$, $\|\Psi_0\|_2=1$ and $\supp \Psi_0\subset B(0,1/R)$, where $R>0$ will be chosen later, we have
$$|V_{\Psi_{(x',\o)}}\sigma(z)|\leq (|V_{\Psi_0} \sigma |\ast|V_{\Psi_{(x',\o)}}\Psi_0|)(z),\quad z\in\bR^{4d},
$$
so that
\begin{align*}
&\|V_{\Psi_{(x',\o)}}\sigma(z(x',\o',x,\o))
\|_{L^\infty_{\o}L^1_{\o'}L^\infty_{x'}L^1_{x}}\\
&\leq \int \|V_{\Psi_{0}}\sigma(z(x',\o',x,\o)-w)
\|_{L^\infty_{\o}L^1_{\o'}L^\infty_{x'}L^1_{x}}
\sup_{(x',\o)\in\R^{2d}}|V_{\Psi_{(x',\o)}}\Psi_0(w)|dw\\
&\leq\sup_{w\in\R^{4d}} \|V_{\Psi_{0}}\sigma(z(x',\o',x,\o)-w)
\|_{L^\infty_{\o}L^1_{\o'}L^\infty_{x'}L^1_{x}} \int
\sup_{(x',\o)\in\R^{2d}}|V_{\Psi_{(x',\o)}} \Psi_0(w)|dw.
\end{align*}
In view of Lemma
\ref{changewind2} we are
therefore reduced to proving
the estimate
\[
\|V_{\Psi_{0}}\sigma(z(x',\o',x,\o)-w)
\|_{L^\infty_{\o}L^1_{\o'}L^\infty_{x'}L^1_{x}}\leq
C\|\sigma\|_{M^{\infty,1}},
\]
uniformly with respect to
$w\in\R^{4d}$. Since,
 \begin{equation*}
 |V_{\Psi_{0}}\sigma(z(x',\o',x,\o)-w)
|\leq \sup_{u\in\R^{2d}}|V_{\Psi_{(x',\o)}}\sigma(u,\tilde{z}
(x',\o',x,\o,w_2))|,
 \end{equation*}
with
 $$\tilde{z}
(x',\o',x,\o,w_2)=(\o'-\nabla_x\Phi(x',\o),x-\nabla_\eta
\Phi(x',\o))-w_2 ,\quad
 w=(w_1,w_2)\in\bR^{4d},
 $$
we shall  prove that
 \begin{equation}\label{T2e}
 \|\sup_{u\in\R^{2d}}|V_{\Psi_{0}}\sigma(u,\tilde{z}
(x',\o',x,\o,w_2)|\|_{L^\infty_{\o}L^1_{\o'}L^\infty_{x'}L^1_{x}}\leq
C\|\sigma\|_{M^{\infty,1}},
 \end{equation}
 uniformly with respect to
  $w_2\in \bR^{2d}$. A
  translation shows that the
  left-hand side in \eqref{T2e}
  is indeed independent of
  $w_2$, and coincides with
  \[
  \sup_{\o\in\R^d}\int_{\R^d_{\o'}}\sup_{x'\in\R^d}\int_{\R^d_x}
\sup_{u\in\R^{2d}} |V_{\Psi_0}\sigma(u;\o'-\nabla_x\Phi(x',\o),x-
\nabla_\eta \Phi(x',\o))|dx\,d\o'.
\]
Here we perform the change of variables $x\longmapsto
x-\nabla_\eta\phi(x',\omega)$. The last expression will be
\[
\leq\sup_{\o\in\R^d}\int_{\R^d_{\o'}}\int_{\R^d_x}
\sup_{u\in\R^{2d}}\sup_{x'\in\R^d}
|V_{\Psi_0}\sigma(u;\o'-\nabla_x\Phi(x',\o),x)|dx\,d\o'.
\]
Now we observe that
\[
\nabla_x\Phi(x',\o)=\nabla_x\Phi(0,\o)+A(x',\o),
\]
where, by the assumption \eqref{fase}, $|A(x',\o)|\leq R$ for some
$R>0$. By Proposition \ref{osservazione3} we can continue the
majorization as
\[
\lesssim \sup_{\o\in\R^d}\int_{\R^d_{\o'}}\int_{\R^d_x}
\sup_{u\in\R^{2d}}\int_{R^d_v}
\varphi_{\o'-\nabla_x\Phi(0,\o),R}(v)
 |V_{\Psi_0}\sigma(u;v,x)|dv\,dx\,d\o'.
\]
Now we bring the supremum
with respect to $u$ inside
the more interior integral
and perform the change of
variables $\o'\longmapsto
\o'-\nabla_x\phi(0,\o)$,
obtaining
\[
\leq \int_{\R^d_{\o'}}\int_{\R^d_x} \int_{R^d_v}
\varphi_{\o',R}(v)
 \sup_{u\in\R^{2d}}|V_{\Psi_0}\sigma(u;v,x)|dv\,dx\,d\o'.
\]
Finally we exchange the two
more interior integral and
apply Lemma
\ref{osservazione2}. The last
expression is seen to be
\[
=C_R\int_{\R^d_{\o'}}\int_{\R^d_x}
 \sup_{u\in\R^{2d}}|V_{\Psi_0}\sigma(u;\o',x)|dx\,d\o'=C_R\|\sigma\|_{M^{\infty,1}}.
 \]
We now prove \eqref{schur2}.
Here we will use repeatedly
the Remark
\ref{cambiovariabili}. By
arguing as above we are
reduced to proving that
\begin{multline}\label{pl}
\sup_{\o'\in\R^{d}}\int_{\R^d_\o}\sup_{x\in\R^d}
\int_{R^d_{x'}}\sup_{u\in\R^{2d}}
|V_{\Psi_0}\sigma(u;\o'-\nabla_x\Phi(x',\o),x- \nabla_\eta
\Phi(x',\o))|dx'\,d\omega\\
\leq
C\|\sigma\|_{M^{\infty,1}}.
\end{multline}
By performing the change of variables
$x'\longmapsto\nabla_\eta\Phi(x',\o)$ and a subsequent
translation, the left-hand side in \eqref{pl} is seen to be
\[
\lesssim \sup_{\o'\in\R^{d}}\int_{\R^d_\o}\sup_{x\in\R^d}
\int_{R^d_{x'}}\sup_{u\in\R^{2d}}
|V_{\Psi_0}\sigma(u;\o'-\nabla_x\Phi(B(x'+x,\o),\o),x')|dx'\,d\omega,
\]
for a suitable function $B$ coming from the inverse change of
variables. Now, by the assumption \eqref{fase} we have
$\nabla_x\Phi(B(x'+x,\o),\o)=\nabla_x\Phi(0,\o)+A'(x,x',\o)$, with
$|A'(x,x',\o)|\leq R$ for a suitable $R>0$. It turns out that the
last expression is
\[
\leq \sup_{\o'\in\R^{d}}\int_{\R^d_\o}
\int_{R^d_{x'}}\sup_{u\in\R^{2d}}\sup_{v\in B(0,R)} |V_{\Psi_0}
\sigma(u;\o'-\nabla_x\Phi(0,\o)+v,x')|dx'\,d\omega.
\]
By the Proposition
\ref{osservazione3} this is
\[
\lesssim \sup_{\o'\in\R^{d}}\int_{\R^d_\o}
\int_{R^d_{x'}}\sup_{u\in\R^{2d}}\int_{\R^d_v}
\varphi_{\o'-\nabla_x\Phi(0,\o),R}(v) |V_{\Psi_0}
\sigma(u;v,x')|dv\,dx'\,d\omega.
\]
Now we perform the change of
variable
$\o\longmapsto\nabla_x\Phi(0,\o)$,
and a subsequent translation,
obtaining
\[
\lesssim \int_{\R^d_\o}
\int_{R^d_{x'}}\sup_{u\in\R^{2d}}\int_{\R^d_v} \varphi_{\o,R}(v)
|V_{\Psi_0} \sigma(u;v,x')|dv\,dx'\,d\omega.
\]
Finally we can bring the
supremum with respect to $u$
inside the more interior
integral, exchange the
integrals with respect to $v$
and $\o'$ and apply Lemma
\ref{osservazione2},
obtaining
\[
\leq C_R \int_{\R^d_\o} \int_{R^d_{x'}}\sup_{u\in\R^{2d}}
|V_{\Psi_0}
\sigma(u;\o,x')|dx'\,d\omega=C_R\|\sigma\|_{M^{\infty,1}},
\]
as desired
\end{proof}\\
We now prove Theorem
\ref{main}. To this end, we
first consider the cases
$(p,q)=(\infty,1)$ (Theorem
\ref{cont2}) and
$(p,q)=(1,\infty)$ (Theorem
\ref{cont2bis}).
\begin{theorem}\label{cont2}
Consider a phase $\Phi$
satisfying {\rm (i)}, {\rm
(ii)} and {\rm (iii)}, and a
symbol $\sigma\in
M_{v_{s,0}\otimes1}^{\infty,1}(\rdd)$,
with $s>d$. Then the
corresponding FIO $T$ extends
to a bounded operator
$\mathcal{M}^{\infty,1}(\rd)\to\mathcal{M}^{\infty,1}(\rd)$.
\end{theorem}
\begin{proof}[Proof] By arguing as at the
beginning of the proof of Theorem
\ref{cont}, we see that it suffices to
prove the estimate
\[
  \sup_{\o\in\R^d}\int_{\R^d_{\o'}}
  \sup_{x'\in\R^d}\int_{\R^d_x}
|V_{\Psi_0}\sigma(x',\o;\o'-\nabla_x\Phi(x',\o),x-
\nabla_\eta \Phi(x',\o))|dx\,d\o'\leq
C\|\sigma\|_{M^{\infty,1}_{v_{s,0}\otimes1}}.
\]
The left-hand side of this estimate is
seen to be
\begin{multline*}
\leq\sup_{\o\in\R^d}\int_{\R^d_{\o'}}
  \sup_{x'\in\R^d}\int_{\R^d_x}
\sup_{(u_1,u_2)\in\R^{2d}}\langle
x'\rangle^{-s}\\
\times
|\langle{u_1}\rangle^{s}V_{\Psi_0}\sigma(u_1,u_2;\o'-\nabla_x\Phi(x',\o),x-
\nabla_\eta \Phi(x',\o))|dx\,d\o',
\end{multline*}
which coincides with
\begin{align*}
&\sup_{\o\in\R^d}\int_{\R^d_{\o'}}
  \sup_{x'\in\R^d}\int_{\R^d_x}
\sup_{(u_1,u_2)\in\R^{2d}}\langle
x'\rangle^{-s}
|\langle{u_1}\rangle^{s}V_{\Psi_0}
\sigma(u_1,u_2;\o'-\nabla_x\Phi(x',\o),x)|dx\,d\o'\\
&\leq \sup_{\o\in\R^d}\int_{\R^d_{\o'}}
  \int_{\R^d_x}
\sup_{(u_1,u_2)\in\R^{2d}}
\langle{u_1}\rangle^{s}
\sup_{x'\in\R^d}\langle x'\rangle^{-s}
|V_{\Psi_0}
\sigma(u_1,u_2;\o'-\nabla_x\Phi(x',\o),x)|dx\,d\o'
\end{align*}
Now we apply Lemma \ref{lemma} to the
function
\[
f(\zeta):=V_{\Psi_0}
\sigma(u_1,u_2;\o'-\nabla_x\Phi(0,\o)-\zeta,x)
\]
with
$v(x')=\nabla_x\Phi(x',\o)-\nabla_x\Phi(0,\o)$
(so that $f(v(x'))=V_{\Psi_0}
\sigma(u_1,u_2;\o'-\nabla_x\Phi(x',\o),x)$).
The assumptions are indeed
satisfied uniformly with
respect to all parameters; in
particular
\[
|v(x')|=|\nabla_x\Phi(x',\o)-\nabla_x\Phi(0,\o)|\leq
C|x'|
\]
for every $(x',\o)\in\R^{2d}$ by
(ii).\\
It follows that we can continue our
majorization as
\begin{align*}
\lesssim
\sup_{\o\in\R^d}\int_{\R^d_{\o'}}
  \int_{\R^d_x}
\sup_{(u_1,u_2)\in\R^{2d}}
\langle{u_1}\rangle^{s}
\int_{\R^d_{x'}}\langle x'\rangle^{-s}
|V_{\Psi_0}
\sigma(u_1,u_2;\o'-\nabla_x\Phi(0,\o)-x',x)|dx'\,dx\,d\o'\\
\leq \sup_{\o\in\R^d}
 \int_{\R^d_{x'}} \langle
x'\rangle^{-s}\int_{\R^d_{\o'}}\int_{\R^d_x}
\sup_{(u_1,u_2)\in\R^{2d}}
\langle{u_1}\rangle^{s}  |V_{\Psi_0}
\sigma(u_1,u_2;\o'-\nabla_x\Phi(0,\o)-x',x)|dx\,d\o'\,dx'.
\end{align*}
By performing the translation
$\o'\longmapsto
\o'-\nabla_x\Phi(0,\o)-x'$ and using
the integrability condition $s>d$ one
sees that this last expression is
\[
=C
 \int_{\R^d_{\o'}}\int_{\R^d_x}
\sup_{(u_1,u_2)\in\R^{2d}}
\langle{u_1}\rangle^{s}
|V_{\Psi_0}
\sigma(u_1,u_2;\o',x)|dx\,d\o'=C\|\sigma
\|_{M^{\infty,1}_{v_{s,0}\otimes1}}.
\]
This concludes the proof.
\end{proof}
\begin{theorem}\label{cont2bis}
Consider a phase $\Phi$
satisfying {\rm (i)}, {\rm
(ii)}, and {\rm (iii)} and a
symbol $\sigma\in
M_{v_{0,s}\otimes1}^{\infty,1}(\rdd)$,
 with $s>d$.
Then, the corresponding FIO
$T$ extends to a bounded
operator
$\mathcal{M}^{1,\infty}(\rd)\to\mathcal{M}^{1,\infty}(\rd)$.
\end{theorem}
\begin{proof}
By arguing as at the beginning of the proof of Theorem \ref{cont2},
we see that it suffices to prove the estimate
\[
  \sup_{\o'\in\R^d}\int_{\R^d_{\o}}
  \sup_{x\in\R^d}\int_{\R^d_{x'}}
|V_{\Psi_0}\sigma(x',\o;\o'-\nabla_x\Phi(x',\o),x- \nabla_\eta
\Phi(x',\o))|dx'\,d\o\leq C\|\sigma\|_{M^{\infty,1}_{v_{0,s}\otimes1}}.
\]
The left-hand side of this estimate is seen to be
\begin{multline*}
\leq\sup_{\o'\in\R^d}\int_{\R^d_{\o}}
  \sup_{x\in\R^d}\int_{\R^d_ {x'}}
\sup_{(u_1,u_2)\in\R^{2d}}\langle
\o\rangle^{-s}\\
\times
|\langle{u_2}\rangle^{s}V_{\Psi_0}\sigma(u_1,u_2;\o'-\nabla_x\Phi(x',\o),x-
\nabla_\eta
\Phi(x',\o))|dx'\,d\o.
\end{multline*}
By performing the change of variables
$x'\longmapsto\nabla_\eta\Phi(x',\o)$ (see Remark \ref{cambiovariabili}) and a subsequent
translation, we obtain
\begin{multline*}
\leq\sup_{\o'\in\R^d}\int_{\R^d_{\o}}
  \sup_{x\in\R^d}\int_{\R^d_ {x'}}
\sup_{(u_1,u_2)\in\R^{2d}}\langle
\o\rangle^{-s}\\
\times
|\langle{u_2}\rangle^{s}V_{\Psi_0}\sigma(u_1,u_2;\o'-\nabla_x\Phi(B(x+x',\o),\o),x')|dx'\,d\o,
\end{multline*}
for a suitable function $B(x',\o)$ coming from the inverse change
of variables. Bringing the supremum over $\o'$ inside,
\begin{multline}\label{mas}
\leq\int_{\R^d_{\o}}
  \sup_{x\in\R^d}\int_{\R^d_ {x'}}
\sup_{(u_1,u_2)\in\R^{2d}}\langle
\o\rangle^{-s}\\
\times \langle{u_2}\rangle^{s} \sup_{\o'\in\R^d}
|V_{\Psi_0}\sigma(u_1,u_2;\o'-\nabla_x\Phi(B(x+x',\o),\o),x')|dx'\,d\o.
\end{multline}
Using Bernstein's inequality \eqref{bernstein}:
$$\sup_{\o'\in\R^d} |f(\o')|\lesssim\intrd |f(\o')|\,d\o',
$$
for
$f(\o')=V_{\Psi_0}\sigma(u_1,u_2;\o'-\nabla_x\Phi(B(x+x',\o),\o),x')$,
and the translation invariance of the Lebesgue measure $d\o'$, the
expression \eqref{mas} is less than
\begin{equation*}
\int_{\R^d_{\o}}\langle \o\rangle^{-s}d\o
  \int_{\R^d_ {x'}}\int_{\R^d_ {\o'}}
\sup_{(u_1,u_2)\in\R^{2d}}\langle{u_2}\rangle^{s}
|V_{\Psi_0}\sigma(u_1,u_2;\o',x')|dx'\,d\o'.
\end{equation*}
Using the integrability condition $s>d$, this last
expression is equal to
\[
C
 \int_{\R^d_ {x'}}\int_{\R^d_ {\o'}}
\sup_{(u_1,u_2)\in\R^{2d}}
\langle{u_2}\rangle^{s}
|V_{\Psi_0}\sigma(u_1,u_2;\o',x')|dx'\,d\o=C\|\sigma
\|_{M_{v_{0,s}\otimes1}^{\infty,1}},
\]
as desired.
\end{proof}

\begin{proof}[Proof of Theorem \ref{main}] Item $(i)$ was already proved in \cite[Theorem 6.1]{fio1}, see also \cite[Theorem 2.1]{toft07}.

Items $(ii)-(iii)$.  By Proposition \ref{proposizione1},  the
  conclusion in Theorem
  \ref{main} is equivalent to
  saying that any FIO
  $\tilde{T}$, with phase
  $\Phi$ satisfying (i),
  (ii) and (iii) and symbol
  $\tilde{\sigma}\in
  M^{\infty,1}(\rdd)$, extends to a
  bounded operator
  $\cM^{p,q}_{v_{0,s_2}}(\rd)\to
  \cM^{p,q}_{v_{-s_1,0}}(\rd)$, for
  $s_1$, $s_2$ as in the
  statement, and moreover
  \[
\|\tilde{T}f\|_{\cM^{p,q}_{v_{-s_1,0}}}\leq
C\|\tilde{\sigma}\|_{M^{\infty,1}}\|f\|_{\cM^{p,q}_{v_{0,s_2}}}.
  \]
  Since this was already
  proved for
  $(p,q)=(\infty,1)$, for all
  $1\leq p=q\leq\infty$, and
  for $(p,q)=(1,\infty)$, the
  desired result follows from
  Proposition \ref{interp}.\end{proof}

Theorem \ref{main} has the
following counterpart in the
framework of Wiener amalgam
spaces. Of course, we are
interested in
$\mathcal{W}(\cF L^p,L^q)$,
with $p\not=q$, since
$\mathcal{W}(\cF
L^p,L^p)=\mathcal{M}^{p}$.

\begin{corollary}\label{mainwiener}
Consider a phase $\Phi$ and a
symbol $\sigma$ as in Theorem
\ref{main}.\par\noindent
$(i)$ Let $1\leq
q<p\leq\infty$. If
$s_1>d\left(\frac{1}{q}-\frac{1}{p}\right)$,
$s_2\geq0$,
 $T$ extends to a
bounded operator on
$\mathcal{W}(\cF L^p,L^q)(\rd)$.
\par\noindent $(ii)$  Let
$1\leq p<q\leq\infty$. If
$s_1\geq0$,
$s_2>d\left(\frac{1}{p}-\frac{1}{q}\right)$,
$T$ extends to a bounded
operator on $\mathcal{W}(\cF
L^p,L^q)(\rd)$.\par
 In all cases,
\begin{equation}\label{normaoperatore2}
\|Tf\|_{\mathcal{W}(\cF L^p,L^q)}\lesssim\|\sigma\|_{M_{v_{s_1,s_2}\otimes1}^{\infty,1}}\|f\|_{\mathcal{W}(\cF L^p,L^q)}.
\end{equation}
\end{corollary}
\begin{proof}
The result easily follows
from Theorem \ref{main}.
Indeed, consider an operator
$T$, with symbol $\sigma$ and
phase $\Phi$, satisfying the
assumptions of Corollary
\ref{mainwiener}. Conjugating
with the Fourier transform
yields the operator
\[
\tilde{T}f(x)=\Fur\circ
T\circ \Fur^{-1}f(x)
\]
Since $\mathcal{M}^{p,q}=\cF^{-1}\mathcal{W}(\cF L^p, L^q),$ it
suffices to prove that $\tilde{T}$ extends to a bounded operator
on $\mathcal{M}^{p,q}$. By duality and an explicit computation
this is equivalent to verifying that the operator
\[
\tilde{T}^\ast f(x)=\intrd
e^{-2\pi
i\Phi(\eta,x)}\overline{\sigma(\eta,x)}\hat{f}(\eta)\,
d\eta
\] extends to a bounded operator
on $\mathcal{M}^{p',q'}$.
Since $\sigma\in
M^{\infty,1}_{v_{s_1,s_2}\otimes
1}(\rdd)$ implies that
$\sigma^\ast(x,\eta)=\overline{\sigma(\eta,x)}\in
M^{\infty,1}_{v_{s_2,s_1}\otimes
1}(\rdd)$ by Lemma \ref{lemm3}, the
desired result is attained
from Theorem \ref{main}.
\end{proof}
\section{Boundedness of
FIOs on  weighted $\mathcal{M}^{p,q}$}
 \label{continuitasumpqw}
Thanks to the commutativity of the
diagram \eqref{diagramma},  the results
of Theorem \ref{main} may be
equivalently stated  as the action of a
FIO $T$ on weighted modulation spaces.
Namely, we have the following result.
\begin{theorem}\label{mainw}
Consider a phase $\Phi$
satisfying {\rm (i)}, {\rm
(ii)} and {\rm (iii)}, and a
symbol $\sigma\in
M^{\infty,1}(\rdd)$. \\
\noindent
$(i)$ If $1\leq
q<p\leq\infty$,  $s_1<-d\left(\frac{1}{q}-\frac{1}{p}\right)$, and $s_2\geq0$, then  $T$ extends to a
bounded operator from $\mathcal{M}^{p,q}_{v_{0,s_2}}(\rd)$ to $\mathcal{M}^{p,q}_{v_{s_1,0}}(\rd)$.
\par\noindent $(ii)$  If
$1\leq p<q\leq\infty$, $s_1\leq0$, and $s_2>d\left(\frac{1}{p}-\frac{1}{q}\right)$,
$T$ extends to a bounded
operator from $\mathcal{M}^{p,q}_{v_{0,s_2}}(\rd)$ to $\mathcal{M}^{p,q}_{v_{s_1,0}}(\rd)$.\par
 In both cases,
\begin{equation*}
\|Tf\|_{\mathcal{M}^{p,q}_{v_{s_1,0}}}\lesssim\|\sigma\|_{M^{\infty,1}}\|f\|_{\mathcal{M}^{p,q}_{v_{0,s_2}}}.
\end{equation*}
\end{theorem}
One can easily rephrase the results of Corollary \ref{mainwiener}
in term of weighted Wiener amalgam spaces.

 We now study the weighted cases not contained above.
\begin{theorem}\label{cont2g}
Consider a phase $\Phi$ satisfying {\rm (i)}, {\rm (ii)} and {\rm
(iii)}, and a symbol $\sigma\in M^{\infty,1}(\rdd)$.  Then, the
corresponding FIO $T$ extends to a bounded operator between the
following spaces:
$$\mathcal{M}^{\infty,1}(\rd)\to\mathcal{M}^{\infty,1}_{v_{0,s}}(\rd),\quad\quad\quad \mathcal{M}^{1,\infty}(\rd)\to\mathcal{M}^{1,\infty}_{v_{s,0}}(\rd),$$
 with $s<-d$, and its norm is bounded from above by $C\|\sigma\|_{M^{\infty,1}}$, for a suitable $C>0$.
\end{theorem}
\begin{proof} To  prove the boundedness of $T$ from $\mathcal{M}^{\infty,1}(\rd)$ to $\mathcal{M}^{\infty,1}_{v_{0,s}}(\rd)$, we have to show that the integral kernel
\[
K_T(x',\omega';x,\omega)=\langle
Tg_{x,\omega},
g_{x',\o'}\ra \la\o'\ra^{s}
\]
satisfies $K_T\in
L^\infty_{\o}L^1_{\o'}L^\infty_{x'}L^1_{x}$. The arguments are similar to those of Theorem \ref{cont2}, we sketch them for sake of clarity. The quantity
\[
  \sup_{\o\in\R^d}\int_{\R^d_{\o'}}
  \sup_{x'\in\R^d}\int_{\R^d_x}
|V_{\Psi_0}\sigma(x',\o;\o'-\nabla_x\Phi(x',\o),x-
\nabla_\eta \Phi(x',\o))|\la \o'\ra^s dx d\o'
\]
can be controlled from above by
\begin{align*}
& C\sup_{\o\in\R^d}\int_{\R^d_{\o'}}
  \sup_{x'\in\R^d}\int_{\R^d_x}
\sup_{(u_1,u_2)\in\R^{2d}} |V_{\Psi_0}\sigma(u_1,u_2;\o'-\nabla_x\Phi(x',\o),x)|\la \o'\ra^sdx\,d\o',\\
&\leq C \sup_{\o\in\R^d}\int_{\R^d_{\o'}}
  \sup_{x'\in\R^d}\int_{\R^d_x}
\sup_{(u_1,u_2)\in\R^{2d}}\sup_{\zeta\in\R^d}|V_{\Psi_0}\sigma(u_1,u_2;\zeta,x)|\la \o'\ra^sdx\,d\o'\\
&\leq C' \sup_{\o\in\R^d}\int_{\R^d_{\o'}}
  \sup_{x'\in\R^d}\int_{\R^d_x}
\sup_{(u_1,u_2)\in\R^{2d}}\int_{\R^d_{\zeta}} |V_{\Psi_0}\sigma(u_1,u_2;\zeta,x)|d\zeta \la \o'\ra^sdx\,d\o',
\end{align*}
where in the last estimate we applied Bernstein's inequality \eqref{bernstein} to the function $f(\zeta)=V_{\Psi_0}\sigma(u_1,u_2;\zeta,x)$. The last estimate, since $s<-d$, is dominated from above by $\|\sigma\|_{M^{\infty,1}}$, as desired.

Similarly, to show $T$: $\mathcal{M}^{1,\infty}\to\mathcal{M}^{1,\infty}_{v_{s,0}}$,  is equivalent to proving that
the integral kernel
\[
K_T(x',\omega';x,\omega)=\langle
Tg_{x,\omega},
g_{x',\o'}\ra \la x'\ra^{s}
\]
satisfies $K_T\in L^\infty_{\o'}L^1_{\o} L^\infty_{x}L^1_{x'}$. The arguments are quite similar to those of Theorem
\ref{cont2bis}. Again, we sketch the proof for sake of clarity.
By performing the change of variables
$\o\longmapsto\nabla_x\Phi(x',\o)$, the quantity
\[
  \sup_{\o'\in\R^d}
  \int_{\R^d_{x'}}\int_{\R^d_{\o}}\sup_{x\in\R^d}
|V_{\Psi_0}\sigma(x',\o;\o'-\nabla_x\Phi(x',\o),x- \nabla_\eta
\Phi(x',\o))| \la x'\ra^s d\o   \,dx'
\]
is less than
\begin{equation*}
\sup_{\o'\in\R^d}\int_{\R^d_ {x'}}\int_{\R^d_{\o}}
\sup_{x\in\R^d}\sup_{(u_1,u_2)\in\R^{2d}}|V_{\Psi_0}\sigma(u_1,u_2;\o,x-\nabla_\eta
\Phi(x',B(x',\o+\o'))|\la x'\ra^s dx'\,d\o,
\end{equation*}
for a suitable function $B(x',\o)$ coming from the inverse change
of variables. The result is attained using
Bernstein's inequality \eqref{bernstein}: $\sup_{x\in\R^d}
|f(x)|\lesssim\intrd |f(x)|\,dx,$ for
$f(x)=V_{\Psi_0}\sigma(u_1,u_2;\o,x)$, and the condition $s<-d$.
\end{proof}

Since the boundedness of the FIO $T$ on $\mathcal{M}^{p}$ is provided
by  Theorem \ref{main}, the complex interpolation with
the preceding results yields:
\begin{theorem}\label{mainw2}
Consider a phase $\Phi$ satisfying {\rm (i)}, {\rm (ii)} and {\rm
(iii)},  a symbol $\sigma\in M^{\infty,1}(\rd)$. \\
\noindent $(i)$ If $1\leq q<p\leq\infty$,  $s_1\geq0$, and
$s_2<-d\left(\frac{1}{q}-\frac{1}{p}\right)$, then $T$ extends to
a bounded operator from $\mathcal{M}^{p,q}_{v_{s_1,0}}(\rd)$ to
$\mathcal{M}^{p,q}_{v_{0,s_2}}(\rd)$.\\
\noindent  $(ii)$  If $1\leq p<q\leq\infty$,
$s_1<-d\left(\frac{1}{p}-\frac{1}{q}\right)$, and $s_2\geq 0$,
then $T$ extends to a bounded operator from
$\mathcal{M}^{p,q}_{v_{0,s_2}}(\rd)$ to
$\mathcal{M}^{p,q}_{v_{s_1,0}}(\rd)$.\par In both cases the norm of $T$
is bounded from above by $C \|\sigma\|_{M^{\infty,1}}$, for a
suitable $C>0$.
\end{theorem}

The results for Wiener amalgam spaces are obtained by similar
arguments as those in Corollary \ref{mainwiener} and left to the
reader.

\section{Boundedness of FIOs
$\mathcal{M}^{p,q}\to
\mathcal{M}^{q,p}$,\ $p\geq q$}
\label{continuitasumqp} In this section
we shall prove the boundedness of an
operator $T$ between
$\mathcal{M}^{p,q}\to
\mathcal{M}^{q,p}$, namely Theorem
\ref{corollario}. As a byproduct,
conditions for the boundedness from
$\mathcal{W}(\cF L^p,L^q)$ to
$\mathcal{W}(\cF L^q,L^p)$ are attained
as well (Corollary \ref{corgen}). By
using complex interpolation, as in the
proof of Theorem \ref{main}, we see
that it suffices to prove the desired
results for $(p,q)=(\infty,1)$:
\begin{theorem}\label{minfu}
 Consider a phase
$\Phi(x,\eta)$ satisfying
{\rm (i)} and  {\rm (ii)}.
Moreover, assume one of the following conditions:\par
(a) the symbol $\sigma\in
M^{\infty,1}(\rdd)$ and, for some
$\delta>0$,
\begin{equation}\label{detcondbis}
   \left|{\rm det}\,
\left(\frac{\partial^2\Phi}{\partial x_i\partial x_l}\Big|_{
(x,\eta)}\right)\right|\geq \delta\quad \forall (x,\eta)\in
\R^{2d},
\end{equation}\par (b) the symbol
$\sigma\in
M^{\infty,1}_{v_{s,0}\otimes1}(\rdd)$,
with $s>d$.\par\noindent
  Then, the
corresponding FIO $T$ extends
to a bounded operator
$\mathcal{M}^{\infty,1}(\rd)\to
\mathcal{M}^{1,\infty}(\rd)$.\par
\end{theorem}
\begin{proof}
We argue as in the first part of the proof of Theorem \ref{cont}.
We see that it suffices to prove that the map $K_T$ defined by
\[
K_T
G(x',\omega')=\int_{\R^{2d}}\langle
Tg_{x,\omega},
g_{x',\o'}\rangle
G(x,\omega)dx\,d\omega
\]
is continuous $L^1_\omega
L^\infty_x \to
L^\infty_\omega L^1_x$. By
Proposition \ref{proschur} it
suffices to prove that its
integral kernel
\[
K_T(x',\omega';x,\omega)=\langle
Tg_{x,\omega},
g_{x',\o'}\rangle
\]
satisfies
\begin{equation}\label{schur1bis}
K_T\in
L^\infty_{\o,\o'}L^1_{x,x'}.
\end{equation}
$(a)$ The same arguments as in the proof of Theorem \ref{cont}
show that it is enough to verify the estimate
\begin{multline}\label{plbis}
\sup_{(\o,\o')\in\R^{2d}}\iint_{\R^{2d}_{x,x'}}
\sup_{u\in\R^{2d}}
|V_{\Psi_0}\sigma(u;\o'-\nabla_x\Phi(x',\o),x-
\nabla_\eta
\Phi(x',\o))|dx\,dx'\\
\leq
C\|\sigma\|_{M^{\infty,1}}.
\end{multline}
To this end, we first perform
the translation  $x\to x-
\nabla_\eta \Phi(x',\o)$ in
the left-hand side of
\eqref{plbis}, obtaining
\[
\sup_{(\o,\o')\in\R^{2d}}\iint_{\R^{2d}_{x,x'}}
\sup_{u\in\R^{2d}}
|V_{\Psi_0}\sigma(u;\o'-\nabla_x\Phi(x',\o),x)|
dx\,dx'.
\]
Then we perform the change of
variables
$x'\longmapsto\nabla_x\Phi(x',\o)$
(followed by a translation)
which by Remark
\ref{cambiovariabili} is a
global diffeomorphism of
$\R^d$, and whose inverse
Jacobian determinant is
uniformly bounded with
respect to all variables.
Hence the last expression is
seen to be
\[
\lesssim
\int_{\R^{d}_{x'}}\int_{\R^{d}_{x}}
\sup_{u\in\R^{2d}}
|V_{\Psi_0}\sigma(u;x',x)|
dx\,dx'=\|\sigma\|_{M^{\infty,1}}.
\]
\noindent  $(b)$ Similar
arguments as above show that
the result is attained once
we prove the estimate
\begin{multline}\label{plbist}
\sup_{(\o,\o')\in\R^{2d}}\iint_{\R^{2d}_{x,x'}}
|V_{\Psi_0}\sigma(x',\o;\o'-\nabla_x\Phi(x',\o),x-
\nabla_\eta
\Phi(x',\o))|dx\,dx'\\
\leq C\|\sigma\|_{M^{\infty,1}_{v_{s,0}\otimes1}}.
\end{multline}
For $s\geq0$, and performing the change of variables $x\longmapsto
x- \nabla_\eta \Phi(x',\o)$, the left-hand side above is
controlled by
$$\leq \sup_{(\o,\o')\in\R^{2d}}\iint_{\R^{2d}_{x,x'}}\la
x'\ra^{-s}\sup_{(u_1,u_2)\in\R^{2d}}\la u_1\ra^s
|V_{\Psi_0}\sigma(u_1,u_2;\o'-\nabla_x\Phi(x',\o),x)|dx\,dx'
$$
$$\leq \sup_{\o'\in\R^{d}}\iint_{\R^{2d}_{x,x'}}\la
x'\ra^{-s}\sup_{(u_1,u_2)\in\R^{2d}}\la
u_1\ra^s
\sup_{\o\in\R^{d}}|V_{\Psi_0}\sigma
(u_1,u_2;\o'-\nabla_x\Phi(x',\o),x)|dx\,dx'.
$$
Now we apply Lemma \ref{lemma} to the function
\[
f(\zeta):=V_{\Psi_0}
\sigma(u_1,u_2;\o'-\nabla_x\Phi(x',0)-\zeta,x)
\]
with $v(\o)=\nabla_x\Phi(x',\o)-\nabla_x\Phi(x',0)$ (so that
$f(v(\o))=V_{\Psi_0} \sigma(u_1,u_2;\o'-\nabla_x\Phi(x',\o),x)$).
The assumptions are indeed satisfied uniformly with respect to all
parameters; in particular $|v(\o)|\leq C|\o|$, for every
$(x',\o)\in\R^{2d}$ by (ii).\\
Continuing  our majorizations, we obtain
\begin{align*}
\lesssim \sup_{\o'\in\R^d}\int_{\R^d_{\o'}}
  \int_{\R^d_x}
\sup_{(u_1,u_2)\in\R^{2d}} \langle{u_1}\rangle^{s}
\int_{\R^d_{x'}}\langle x'\rangle^{-s} |V_{\Psi_0}
\sigma(u_1,u_2;\o'-\nabla_x\Phi(x',0)-\o,x)|d\o \,dx'\,dx\,\\
\leq \sup_{\o'\in\R^d}
 \int_{\R^d_{x'}} \langle
x'\rangle^{-s}\int_{\R^d_{\o}}\int_{\R^d_x}
\sup_{(u_1,u_2)\in\R^{2d}} \langle{u_1}\rangle^{s}  |V_{\Psi_0}
\sigma(u_1,u_2;\o,x)|\,d\o\,dx\,dx'
\end{align*}
where in the last row we
perform the translation (up
to a sign) $\o\longmapsto
\o'-\nabla_x\Phi(x',0)-\o$ and
using the integrability
condition $s>d$ one sees that
this last expression is
dominated by $\|\sigma
\|_{M^{\infty,1}_{v_{s,0}\otimes1}}.$
\end{proof}\\
 From
Theorem \ref{corollario} and
using the same arguments as
in Corollary
\ref{mainwiener}, one can
prove the following result.
\begin{corollary}\label{corgen}
Let $1\leq q<p\leq\infty$.
Consider a phase
$\Phi(x,\eta)$ satisfying
{\rm (i)}, {\rm (ii)} and
{\rm (iii)}.
 Moreover,
assume one of the following conditions:\par (a)  the symbol $\sigma\in M^{\infty,1}(\rdd)$ and,
for some $\delta>0$,
\begin{equation}\label{detcond01}
   \left|{\rm det}\,
\left(\frac{\partial^2\Phi}{\partial \eta_i\partial \eta_l}\Big|_{
(x,\eta)}\right)\right|\geq \delta\quad \forall (x,\eta)\in
\R^{2d},
\end{equation}
or\par (b) the symbol
$\sigma\in
M^{\infty,1}_{v_{0,s}\otimes1}(\rdd)$,
with
$s>d\left(\frac{1}{q}-\frac{1}{p}\right)$.\par\noindent
  Then, the
corresponding FIO $T$ extends
to a bounded operator
$\mathcal{W}(\Fur L^p,L^q)(\rd)\to
\mathcal{W}(\Fur
L^q,L^p)(\rd)$.\par
\end{corollary}
\begin{example}\label{metaplettici}\rm
(Metaplectic operators)
Consider the particular case
of quadratic phases, namely
phases of the type
\begin{equation}\label{faseq}
\Phi(x,\eta)=\frac{1}{2}
Ax\cdot
x+Bx\cdot\eta+\frac{1}{2}C\eta\cdot\eta+\eta_0\cdot
x-x_0\cdot\eta,
\end{equation}
where
$x_0,\eta_0\in\mathbb{R}^d$,
$A,C$ are real symmetric
$d\times d$ matrices and $B$
is a real $d\times d$
nondegenerate matrix.\par It
is easy to see that, if we
take the symbol
$\sigma\equiv1$ and the phase
\eqref{faseq}, the
 corresponding FIO $T$ is (up to a constant
factor) a metaplectic
operator. This can be seen by
means of the easily verified
factorization
\begin{equation}\label{factorization}
 T=M_{\eta_0} U_A
D_B\Fur^{-1} U_C\Fur T_{x_0},
\end{equation}
where $U_A$ and $U_C$ are the
multiplication operators by
$e^{\pi i Ax\cdot x}$ and
$e^{\pi i C\eta\cdot \eta}$
respectively, and $D_B$ is
the dilation operator
$f\mapsto f(B\cdot)$. Each of
the factors is (up to a
constant factor) a
metaplectic operator (see
e.g. the proof of
\cite[Theorem
18.5.9]{hormander}), so $T$
is. We refer to
\cite[Proposition 2.7
(ii)]{cordero2} and
\cite[Section 7]{fio1} for
details about the symplectic
matrix which yields such an
operator.\par We see that the
assumptions (i) and (ii) in
the Introduction are clearly
satisfied, whereas the
hypotheses (iii) in the
Introduction,
\eqref{detcondbis0} and
\eqref{detcond01} are
equivalent to ${\rm
det}\,B\not=0$, ${\rm
det}\,A\not=0$ and ${\rm
det}\,C\not=0$ respectively.
In particular, the first part
of Corollary \ref{corgen}
generalizes \cite[Theorem
4.1]{cordero2}.
\end{example}

\section{Some counterexamples related to the
Schr\"odinger
multiplier}\label{schroedinger}
 In
this section, we study the
action of the Scr\"odinger
multiplier
$$f\longmapsto\mathcal{F}^{-1}\left(e^{i\pi|\cdot|^2}\hat{f}\right),$$
that is the multiplier with
symbol $e^{i\pi|\eta|^2}$, on
the Wiener amalgam spaces
$W(\Fur L^p,L^q_s)$.
Equivalently, we study the
action of the pointwise
multiplication operator
\begin{equation}\label{moltiplicazione}
Af(x)= e^{i\pi|x|^2}f(x),\quad x\in\rd
\end{equation}
on the modulation spaces
$\mathcal{M}^{p,q}_{v_{0,s}}$.
This will provide useful
examples to test the
sharpness of the thresholds
arising in our results.
Notice that $A$ is the FIO
with phase
$\Phi(x,\eta)=x\eta+\frac{|\eta|^2}{2}$
and symbol $\sigma\equiv1$.
In particular, it satisfies \eqref{detcondbis0}.
 \par It was shown in
\cite[Proposition 7.1]{fio1}
that $A$ is bounded on
$\cM^{p,q}$ (if and) only if
$p=q$ (for a proof of the
positive results see
\cite{benyi,cordero2,fio1}).
It is natural to wonder
whether for suitable negative
values of $s$ it is bounded
as an operator $\cM^{p,q}\to
\cM^{p,q}_{v_{0,s}}$. The
next result deals with the
optimal range of $s$ for this
to happen.
\begin{proposition}\label{3.1}
If the operator $A$ in
\eqref{moltiplicazione} is
bounded as an operator
$\cM^{p,q}(\rd)\to
\cM^{p,q}_{v_{0,s}}(\rd)$, for some
$1\leq p,q\leq\infty$,
$s\in\R$, then $p\geq q$ and
$s\leq -d(1/q-1/p)$.
Moreover, if $p=\infty$ and
$q<\infty$, then $s<-d/q$.
\end{proposition}
\begin{proof}
We will prove later on that $p\geq q$.
Let us verify now the remaining
conditions.
 We
test the estimate
\begin{equation}\label{equazione}
\|Af\|_{M^{p,q}_{v_{0,s}}}\lesssim\|f\|_{M^{p,q}}
\end{equation}
on the family of functions
$f_\lambda(x)=e^{-\pi\lambda|x|^2}$,
$\lambda>0$. Applying Lemma
\ref{lemm2} with $a=\lambda$,
$b=0$ we see that
\begin{equation}\label{secondomembro}
\|f_\lambda\|_{M^{p,q}}\asymp\|G_{a+ib}\|_{W(\Fur
L^p,L^q)}\asymp\lambda^{-\frac{d}
{2p}},\qquad{\rm as}\
\lambda\to0^+.
\end{equation}
On the other hand, with the
notation of Lemma \ref{lemm2},
 we
have
\[
\|Af_\lambda\|_{M^{p,q}_{v_{0,s}}}\asymp\|G_{a+ib}\|_{W(\Fur
L^p,L^q_s)}
\]
with $a=\lambda$, $b=-1$. We now
estimate this last expression. We see
that, by \eqref{lem10},
\[
\|G_{a+ib}T_x g\|_{\Fur
L^p}\asymp
e^{-\frac{\pi\lambda|x|^2}{1+\lambda+\lambda^2}},\qquad
\lambda\leq1.
\]
Let
$\mu=\lambda/(1+\lambda+\lambda^2)$
(observe that
$\mu\sim\lambda$ and $\log
\mu\sim\log\lambda$ as
$\lambda\to0^+$). Then, for
$q<\infty$,
\begin{align}
\|G_{a+ib}\|_{W(\Fur
L^p,L^q_s)}&\asymp\left(\int
e^{-\pi
q\mu|x|^2}(1+|x|)^{qs}\,dx\right)^{1/q}\nonumber\\
&=\mu^{-\frac{d}{2q}}\left(\int
e^{-\pi
q|x|^2}(1+\mu^{-\frac{1}{2}}|x|)^{qs}\,dx\right)^{1/q}.\label{lmn}
\end{align}
First,  assume $s\leq0$. If $\mu^{\frac{1}{2}}\leq|x|\leq1$ we
have $e^{-\pi q|x|^2}\gtrsim 1$ and
$1+\mu^{-\frac{1}{2}}|x|\leq
2\mu^{-\frac{1}{2}}|x|$. Hence the last
expression turns out to be
\begin{align*}
&\geq
\mu^{-\frac{d}{2q}}\left(
\int_{\mu^{\frac{1}{2}}\leq|x|\leq1}
e^{-\pi
q|x|^2}(1+\mu^{-\frac{1}{2}}|x|)^{qs}\,dx\right)^{1/q}\\
&\gtrsim
\mu^{-\frac{d}{2q}-\frac{s}{2}}
\left(\int_{\mu^{\frac{1}{2}}\leq|x|\leq1}
|x|^{qs}\,dx\right)^{1/q}.
\end{align*}
On the other hand, using
polar coordinates one sees
that, as $\mu\to0^+$,
\[
\int_{\mu^{\frac{1}{2}}\leq|x|\leq1}
|x|^{qs}\,dx\asymp
\begin{cases}
1 &{\rm if}\ \ sq>-d\\
|\log\mu|&{\rm if}\ \ sq=-d\\
\mu^{\frac{d}{2}+\frac{sq}{2}}&{\rm
if}\ \ sq<-d.
\end{cases}
\]
Secondly, assume $s>0$ in
\eqref{lmn}. Since
$1+\mu^{-\frac{1}{2}}|x|\geq
\mu^{-\frac{1}{2}}|x|$, we
have
$$\|G_{a+ib}\|_{W(\Fur
L^p,L^q_s)}\gtrsim \mu^{-\frac{d}{2q}-\frac{s}{2}}.$$
Hereby, we deduce, as
$\lambda\to0^+$,
\begin{equation}\label{ala9}
\|Af_\lambda\|_{M^{p,q}_{v_{0,s}}}\gtrsim
\begin{cases}
\lambda^{-\frac{d}{2q}-\frac{s}{2}} &{\rm if}\ \ sq>-d\\
\lambda^{-\frac{d}{2q}-\frac{s}{2}}
|\log\lambda|^\frac{1}{q}&{\rm if}\ \ sq=-d\\
1&{\rm if}\ \ sq<-d.
\end{cases}
\end{equation}
By combining this estimate with
\eqref{equazione} and
\eqref{secondomembro} and letting
$\lambda\to0^+$ we deduce the desired
conclusion, when $q<\infty$. An easy
modification of the above arguments
yields also the case $p=q=\infty$.\par
Let us now prove the condition $p\geq
q$. By duality, it suffices to prove
that if the inequality
\begin{equation}\label{ala0}
\|Af\|_{M^{p',q'}}\leq C\|f\|_{M^{p',q'}_{v_{0,s}}},\qquad\forall
f\in\cS(\R^d),
\end{equation}
holds true for some $s\in\R$, then
$p\geq q$. To see this, we test the
estimates on the Schwartz functions
$f_\lambda$, $0<\lambda<1$, already
used before. It follows from
\eqref{ala9}, in particular, that
\begin{equation}\label{ala}
\|Af_\lambda\|_{M^{p',q'}}\gtrsim
\lambda^{-\frac{d}{2q'}},\qquad
0<\lambda\leq1.
\end{equation}
On the other hand, by \eqref{lem10}
with $a=\lambda,b=0$,
\[
\|G_{a+ib} T_x g\|_{\Fur
L^{p'}}\asymp\lambda^{-\frac{d}{2p'}}e^{-\pi\nu|x|^2},
\qquad 0<\lambda\leq1,
\]
where
$\nu=\lambda/(\lambda^2+\lambda)\asymp
1$ for $0<\lambda\leq1$. Hence
\[
\|f_\lambda\|_{M^{p',q'}_{v_{0,s}}}\lesssim
\lambda^{-\frac{d}{2p'}},\qquad 0<\lambda\leq1.
\] This estimate,
combined with \eqref{ala0} and
\eqref{ala} and letting $\lambda\to0^+$
yields $p'\leq q'$, namely $p\geq q$.
\end{proof}\\
Here is a related example.
\begin{example}{\rm Let us show the unboundedness of the operator $A$ between $M^{1,\infty}(\rd)$ and $M^{\infty,q}(\rd)$, for every $q<\infty$.
 Consider the tempered distribution $\delta$, defined
 by $\la\delta,\f\ra=\bar{\f}(0)$, for every
  $\f\in\cS(\rd)$. Then, for a fixed
   non-zero window $g\in\cS(\rd)$,  the STFT
$V_g \delta\phas=\la \delta, M_\o T_x
g\ra={\overline{g(-x)}}\in
L^{1,\infty}(\rdd),$ that is, $\delta
\in M^{1,\infty}(\rd).$ Now,
\begin{equation}\label{test}
V_g (A\delta)\phas=\la \delta,e^{-\pi i |\cdot|^2} M_\o T_xg\ra= \overline{g(-x)}\notin L^{\infty,q}(\rdd),
\end{equation}
for every $q<\infty$. \par This also shows the unboundedness of
the operator $A$ between $M^{\infty}(\rd)$ and $M^\infty_{v_{0,s}}(\rd)$,
for every $s>0$. Indeed, using the inclusion relations
$$M^\infty_{v_{0,s}}\subset M^{\infty,q},\quad
s>\frac1q,$$ we see that, if $A$ were bounded between $M^{\infty}$
and $M^\infty_{v_{0,s}}$ for some $s>0$, then  the operator $A$
would be bounded between $M^{\infty}$ and $M^{\infty,q}$ as well,
and this is false, as shown by \eqref{test}. }
\end{example}
\begin{proposition}\label{3.1bis}
If the operator $A$ in
\eqref{moltiplicazione} is
bounded as an operator
$\cM^{p,q}(\rd)\to
\cM^{p,q}_{v_{s,0}}(\rd)$ for some
$1\leq p,q\leq\infty$,
$s\in\R$, then $s\leq
-d(1/q-1/p)$. Moreover, if
$p=\infty$ and $q<\infty$,
then $s<-d/q$.
\end{proposition}
We need the following
 elementary result
 \begin{lemma}\label{lemmadis}
 Let $a_1,a_2,C_0\in\R$ with
 $0<a_1\leq C_0$, $C_0^{-1}\leq|a_2|\leq
 C_0$. Then
 \[
 \int_{\R^d}
 e^{-a_1|\omega+a_2x|^2}\langle\omega\rangle^sd\omega\gtrsim
 \langle x\rangle^s,
 \]
 where the constant implicit
 in the notation $\gtrsim$
 only depends on $d,C_0$ and $s$.
\begin{proof}
Let us observe that
$e^{-a_1|\omega+a_2x|^2}\geq
e^{-1}$ for $\omega\in
B(-a_2x,a_1^{-1/2})$. Hence
\begin{equation*}
\int_{\R^d}
 e^{-a_1|\omega+a_2x|^2}\langle\omega\rangle^sd\omega
 \gtrsim \int_{B(0,a_1^{-1/2})}
 \langle y-a_2x\rangle^sdy\gtrsim \langle a_2
 x\rangle^s \int_{B(0,a_1^{-1/2})}
 \langle y\rangle^{-|s|}dy,
\end{equation*} which gives the
conclusion by the assumptions
on $a_1$ and $a_2$.
\end{proof}
 \end{lemma}\begin{proof}[Proof of
Proposition \ref{3.1bis}] In
view of what we proved in
Proposition \ref{3.1}, it
suffices to verify that if
the estimate
\[
\|Af\|_{\cM^{p,q}_{v_{s,0}}}\leq
C\|f\|_{\cM^{p,q}},\qquad\forall
f\in\cS(\R^d),
\]
holds true, then
\eqref{equazione} holds as
well, at least for all
$f(x)=f_\lambda(x)=e^{-\pi\lambda|x|^2}$,
$0<\lambda<1$. Hence we are
left to prove that
\[
\|Af_\lambda\|_{\cM^{p,q}_{v_{s,0}}}\gtrsim
\|Af_\lambda\|_{\cM^{p,q}_{v_{0,s}}},
\]
or equivalently that
\[
\|G_{a+ib}\|_{W(\Fur
L^p_s,L^q)}\gtrsim
\|G_{a+ib}\|_{W(\Fur
L^p,L^q_s)},
\]
with $a=\lambda$, $b=-1$ (see
the notation of Lemma
\ref{lemm2}). To this end,
observe that by
\eqref{lem1000} we have, for
$p<\infty$,
\begin{align*}
&\|G_{a+ib}T_x g\|_{\Fur
L^p_s}\\
&\quad=((a+1)^2+b^2)^{-d/4}
e^{-\frac{\pi(a+1)|x|^2}{(a+1)^2+b^2}}\left(\int
 e^{-\frac{p\pi}{(a+1)^2+b^2}[(a^2+b^2+a)|\omega|^2+2bx\omega]}\langle\omega\rangle^{ps}d\omega\right)^{\frac{1}{p}}\nonumber\\
 &\quad=((a+1)^2+b^2)^{-d/4} e^{-\frac{\pi
 a|x|^2}{a(a+1)+b^2}}\left(\int
 e^{-\frac{p\pi}{(a+1)^2+b^2}|
 \sqrt{a(a+1)+b^2}\omega+
 \frac{b}{\sqrt{a(a+1)+b^2}}x|^2}\langle\omega\rangle^{ps}d\omega
 \right)^\frac{1}{p}\nonumber
\end{align*} Since $a=\lambda\leq 1$ and $b=-1$ an application of
Lemma \ref{lemmadis}
 gives
\[
\|G_{a+ib}T_x g\|_{\Fur
L^p_s}\gtrsim e^{-\frac{\pi
 \lambda|x|^2}{\lambda^2+\lambda+1}}\langle x\rangle^{s}
 \asymp \|G_{a+ib}T_x g\|_{\Fur
L^p}\langle x\rangle^{s},
 \]
 where for the last estimate
 we used \eqref{lem10}.\par
 Similarly one treats  the
 case $p=\infty$.
\end{proof}

\begin{remark}\label{aspettato}\rm
By duality, the assumption in
Theorem \ref{3.1} could be
rephrased as the boundedness
of $A$ itself between
$\cM^{p',q'}_{v_{0,-s}}(\rd)\to
\cM^{p',q'}(\rd)$. Hence, as a
consequence of Proposition
\ref{proposizione1} and
Proposition \ref{3.1}, we see
that the threshold for $s_2$
arising in Theorem \ref{main}
$(iii)$ (with $s_1=0$) is
essentially sharp; namely, it
is sharp for $p=1$, whereas
when $p>1$ only the case of
equality remains open.\\
Similarly, Proposition
\ref{3.1bis} shows that the
threshold for $s_1$ arising
in Theorem \ref{main} $(ii)$
(with $s_2=0$) is essentially
sharp.
\end{remark}
The following result shows
that the conclusion in
Theorem \ref{corollario}
generally fails if $q>p$.
\begin{proposition}
\label{controesempio0} If the
operator $A$ in
\eqref{moltiplicazione} is
bounded as an operator
$\cM^{p,q}(\rd)\to \cM^{q,p}(\rd)$ for
some $1\leq p,q\leq\infty$,
then $q\leq p$.
\end{proposition}
\begin{proof}
We test the estimate
\[
\|Af\|_{\cM^{q,p}}\lesssim
\|f\|_{\cM^{p,q}},\qquad
\forall f\in\cS(\R^d),
\]
 on the family of Schwartz functions $f_\lambda(x)=e^{-\pi\lambda|x|^2}$,
$\lambda>0$. Applying Lemma
\ref{lemm2} with $a=\lambda$,
$b=0$ we see that
\begin{equation}\label{secondomembro2}
\|f_\lambda\|_{M^{p,q}}\asymp\|G_{a+ib}\|_{W(\Fur
L^p,L^q)}\asymp\lambda^{\frac{d}
{2q}-\frac{d}{2}}\qquad{\rm
as}\ \lambda\to+\infty.
\end{equation}
On the other hand, with the
notation of Lemma \ref{lemm2}
with $a=\lambda$, $b=1$,
 we
have
\[
\|Af_\lambda\|_{M^{q,p}_s}\asymp\|G_{a+ib}\|_{W(\Fur
L^q,L^p)}\asymp\lambda^{\frac{d}
{2p}-\frac{d}{2}}\qquad{\rm
as}\ \lambda\to+\infty.
\]
Hence it turns out $q\leq p$.
\end{proof}\\
Finally we present a
counterexample related to
Theorem \ref{corollario}.

\begin{proposition}\label{controesempio}
The Schr\"odinger multiplier
$f\longmapsto\mathcal{F}^{-1}\left(e^{i\pi|\cdot|^2}\hat{f}\right)$
is not bounded as an operator
$\cM^{p,q}(\rd)\to\cM^{q,p}(\rd)$ if
$p\not=q$.
\end{proposition}
\begin{proof}
It suffices to prove that the
pointwise multiplication
operator $A$ in
\eqref{moltiplicazione} is
not bounded as an operator
$\mathcal{W}(\Fur
L^p,L^q)(\rd)\to\mathcal{W}(\Fur
L^q,L^p)(\rd)$ if $p\not=q$.\\
 We
test the estimate
\[
\|Af\|_{\mathcal{W}(\Fur
L^q,L^p)}\leq
C\|f\|_{\mathcal{W}(\Fur
L^p,L^q)},\qquad \forall
f\in\cS(\R^d),
\]
on the family of Schwartz
functions
$f_\lambda(x)=e^{-\pi\lambda|x|^2}$,
$\lambda>0$.\\
By applying Lemma \ref{lemm2}
with $a=1/\lambda$ and $b=0$
we obtain
\[
\|f_\lambda\|_{\mathcal{W}(\Fur
L^p,L^q)}=\lambda^{-d/2}\|G_{a+ib}\|_{\mathcal{W}(\Fur
L^p,L^q)}\asymp
\begin{cases}\lambda^{-\frac{d}{2q}}&0<\lambda\leq1,\\
\lambda^{-\frac{d}{2p'}}&\lambda\geq1.
\end{cases}
\]
Similarly we have, with
$a=\frac{\lambda}{\lambda^2+1}$,
$b=\frac{1}{\lambda^2+1}$,
\[
\|Af_\lambda\|_{\mathcal{W}(\Fur
L^q,L^p)}=(a^2+b^2)^{d/4}\|G_{a+ib}\|_{\mathcal{W}(\Fur
L^q,L^p)}\asymp
\begin{cases}\lambda^{-\frac{d}{2p}}&0<\lambda\leq1,\\
\lambda^{-\frac{d}{2q'}}&\lambda\geq1.
\end{cases}
\]
Letting $\lambda\to0^+$ and
$\lambda\to+\infty$ one
deduces a contradiction
unless $p=q$.
\end{proof}
Notice that this multiplier is the FIO with phase
$\Phi(x,\eta)=x\eta+\frac{|\eta|^2}{2}$ and symbol
$\sigma\equiv1$, so that neither \eqref{detcondbis0} nor condition
$(b)$ is  satisfied (whereas \eqref{detcond01} is). Indeed, for
$\sigma\equiv1$, $g\in\cS(\rdd)$,  we have
 $$|V_g \sigma(z,\zeta)|=|\hat{\bar{g}}(\zeta)|\notin L^1(\R_\zeta^{2d},L^\infty_{v_{s,0}}(\R_z^{2d})),\quad\forall s>0.
 $$


\end{document}